\documentclass[11pt,leqno]{amsart}
\usepackage{amsmath,amsthm}
\usepackage{amsfonts,amssymb} 
\usepackage{amsrefs}
\newtheorem{theorem}{Theorem}
\newtheorem{lemma}[theorem]{Lemma}
\newtheorem{proposition}[theorem]{Proposition}

\begin{document}
\title[Dynamics for the Brownian web]{Dynamics for the Brownian web and the erosion flow}
\author{Chris Howitt}
\address{Department of Statistics, University of Warwick, Coventry CV4 7AL, UK}
\email{c.j.howitt@warwick.ac.uk}
\author{Jon Warren}
\email{j.warren@warwick.ac.uk}
\date{}

\begin{abstract}
The Brownian web is a random object that occurs as the scaling limit of an infinite system of coalescing random walks.
Perturbing this system of random walks by, independently at each point in space-time, resampling the random walk increments, leads to 
some natural dynamics. In this paper we consider  the corresponding dynamics for the Brownian web. In particular, pairs of coupled Brownian webs are studied, where the second web is obtained from the first by perturbing according to these dynamics.  A stochastic flow of kernels, which we call the erosion flow, is obtained via a filtering  construction from such coupled Brownian webs, and  the $N$-point  motions of this flow of kernels are identified. 
\end{abstract}
 
\maketitle

\section{Introduction}

Suppose that $\bigl(\xi_{k,n}; (k,n) \in L \bigr)$ is a family of independent random signs with ${\mathbf P}( \xi_{k,n} =1)={\mathbf P}(\xi_{k,n}=-1 \bigr)=\tfrac{1}{2}$, indexed by the points of the  lattice $L= \bigl\{ (k,n) \in {\mathbf Z}^2: k+n \text{ is even }\bigr\}$. By representing  $\xi_{k,n}=1$ as an arrow pointing from $(k,n)$ to $(k+1,n+1)$ and $\xi_{k,n}=-1$ as an arrow pointing from $(k,n)$ to $(k-1,n+1)$ and joining up the arrows, starting from arbitrary points $(k,n) \in L$, we construct an infinite family, ${\mathcal S}$, of coalescing simple random walk paths. 

The scaling limit of this system of random walks, a coalescing system of Brownian motions, was first investigated by Arratia, \cite{arratia}. A detailed 
study was made by T\'{o}th and Werner, \cite{tw}.
More recently Fontes, Isopi, Newman, and Ravishankar, \cite{FINR}, proposed a  framework in which a configuration of  paths with various starting points in space-time is treated as a point in a suitable metric space ${\mathcal H}$. Letting ${\mathcal S}^\epsilon$ denote the family of paths obtained by scaling ${\mathcal S}$ by a factor of $\epsilon$ in time and $\sqrt{\epsilon}$  in space, they prove that as $\epsilon$ tends to zero, then  ${\mathcal S}^\epsilon$ converges in distribution to ${\mathcal W}$, where the latter is an ${\mathcal H}$-valued random variable they named the Brownian web.
For more on the Brownian web, see \cite{MR2165253}, \cite{MR2196970}, and \cite{MR2120239}. A recent  paper by Sun and Swart, \cite{ss}, contains a construction of pairs of Brownian webs that is closely related to the work presented here.

Suppose that each random sign $\xi_{k,n}$ is replaced by a stochastic process $\bigl(\xi_{k,n}(u); u \in {\mathbf R}\bigr)$ evolving as as stationary Markov chain with state space $\bigl\{-1,+1\bigr\}$ and jumping from $-1$ to $+1$, and vice-versa at unit rate. Assume that these processes, as $(k,n)$ varies through $L$, are mutually independent of one another. 
At each instant $ u \in {\mathbf R}$ we may construct from the random signs $  \bigl( \xi_{k,n}(u); (k,n) \in L \bigr)$ a system of coalescing paths. Representing this system as a point ${\mathcal S}(u)$ in the metric space ${\mathcal H}$ we obtain a stationary ${\mathcal H}$-valued Markov process $\bigl( {\mathcal S}(u); u \in {\mathbf R} \bigr)$.  

It is now natural to consider the behaviour of the scaled processes $\bigl({\mathcal S}^\epsilon(u); u \in {\mathbf R}\bigr)$. It is not difficult to guess that
as $\epsilon $ tends down to zero, then for any $u_1 \neq u_2$ the   law of  the pair $\bigl({\mathcal S}^\epsilon(u_1), {\mathcal S}^\epsilon(u_2)\bigr)$ converges to that of  a pair of independent Brownian webs. In fact this is a manifestation of the sensitivity of the coalescing random walks to perturbations, and  correspondingly  the Brownian web is an example of a black noise in Tsirelson's theory of continuous products of probability spaces, see \cite{stflour} and \cite{tsirelson}.    However if we slow down the speed of the perturbations as we scale, then we may obtain a non-trivial limit. In fact,  the law of the  pair    $\bigl({\mathcal S}^\epsilon(\sqrt{\epsilon} u_1), {\mathcal S}^\epsilon(\sqrt{\epsilon} u_2)\bigr)$ converges to that of  a pair  which we call a $\theta$-coupling of Brownian webs. Here $\theta$ is a non-negative real parameter, which is proportional to $|u_1-u_2|$. If we consider a  pair  of $\theta$-coupled webs $\bigl( {\mathcal W}, {\mathcal W}^\prime\bigr)$, and let $\bigl(X(t); t \geq 0 \bigr)$ denote the, almost surely unique, Brownian path that starts from, say, the origin of space-time in ${\mathcal W}$, and $\bigl(X^\prime(t); t \geq 0 \bigr)$ denote the corresponding path taken from the web ${\mathcal W}^\prime$, then the pair $\bigl(X, X^\prime)$ is distributed according to the law of a diffusion in ${\mathbf R}^2$, which we have previously, \cite{hw}, called $\theta$-coupled Brownian motions.  This means $X$ and $X^\prime$ are both Brownian motions relative to some 
common filtration, and,
\begin{equation}
\langle X,X^\prime \rangle (t)=\int_0^t {\mathbf 1}_{( X(s)=X^\prime(s))} ds  \qquad  t\geq 0, 
\end{equation}
\begin{equation}
L^0_t\bigl( X-X^\prime\bigr) = 2\theta \int_0^t {\mathbf 1}_{( X(s)=X^\prime(s))} ds  \qquad t\geq 0.
\end{equation}
Here $L^0_t(X-X^\prime)$ denotes the semimartingale local time accrued by $X-X^\prime$ at zero by time $t$.
A pair of $\theta$-coupled Brownian motions evolve independently of each other when apart, but when they meet there is some interaction, 
often described as stickiness, which results in their  momentarily moving together.

As a consequence of these observations on scaling, it is reasonable to suppose that there exists a stationary  ${\mathcal H}$-valued Markov process $\bigl({\mathcal W}(u); u \in {\mathbf R}\bigr)$ with, for each $u$, the law of ${\mathcal W}(u) $ being that of the Brownian web, and with the law of $\bigl({\mathcal W}(u_1),{\mathcal W}(u_2) \bigr)$ being that of some pair of $\theta$-coupled webs. However we caution the reader that we limit ourselves here to  considering the finite dimensional distributions  of any such process, and leave to another time a more complete study.

A stochastic flow of kernels on a measurable space $\bigl(E, {\mathcal E}\bigr)$ is a doubly indexed family $\bigl(K_{s,t}; s \leq t\bigr)$ of 
random $E\times {\mathcal E}$ transition kernels satisfying the flow property
\begin{equation}
K_{s,u}(x, A) = \int_E K_{s,t}(x, dy ) K_{t,u}(y,A) \qquad \qquad x \in E, A \in {\mathcal E} 
\end{equation} 
almost surely for each $ s \leq t$.  We also postulate  independent and stationary innovations  in that 
$K_{t_1,t_2}, K_{t_2,t_3}, \ldots K_{t_{n-1},t_n}$ are independent  for all choices of $t_1<t_2<\ldots...<t_n$ and 
$K_{s,t}\stackrel{dist}{=} K_{s+h,t+h}$ for all $s<t$ and $h\in {\mathbf R}$.  The general theory of such flows was developed by Le Jan and Raimond in 
\cite{lejan}, see also Tsirelson \cite{tsirelson}.

One possible interpretation for the flow $K$ is as describing  the random evolution of a distribution of mass on $E$. In this case 
$K_{s,t}(x,A)$ represents the proportion of that mass which  was located at $x$ at time $s$ which is within the set $A$ at time $t$. An 
alternative interpretation of $K$ is as a  random environment- in time and space- governing the motion of a particle. Then $K_{s,t}(x,A)$ 
is the conditional probability given the environment that a particle which is located at $x$ at time $s$ is located within the set $A$ at 
time $t$.  

Consider once again the stationary process $\bigl({\mathcal S}(u) \in {\mathbf R}\bigr)$ generated from the random sign processes $ \bigl(\xi_{k,n}(u); u \in {\mathbf R} \bigr)$ for $(k,n) \in L$. Suppose that  $u_1\neq u_2$ and let ${\mathcal S}={\mathcal S}(u_1)$ and ${\mathcal S}^\prime={\mathcal S}(u_2)$. From the pair $\bigl({\mathcal S},{\mathcal S}^\prime\bigr)$ can naturally define a stochastic flow of kernels on the integer lattice  via
\begin{equation}
K_{m,n}(k,A)= {\mathbf P} \bigl ( S^\prime_{k,m}(n) \in A | {\mathcal S} \bigr)
\end{equation}
for  $A \subseteq {\mathbf Z}$ and $(k,m) \in L$, where  $S^\prime_{k,m}(n)$ denotes the position at time $n\geq m$ of the path in ${\mathcal S}^\prime$ starting from $(k,m)$.
Actually this  doesnt quite fit the definition of a flow given above because of the periodic nature of the simple random walk, but this is not of any 
importance.  
The corresponding flow of mass on the integer lattice ${\mathbf Z}$ is easily described.  Whatever mass is at $k$ at time $n$ split into two, with a proportion $p$ of it moving to $k+\xi_{k,n}$ at time $(n+1)$, whilst the reminder moves to $k-\xi_{k,n}$.  Here $p={\mathbf P}\bigl( \xi_{k,n}^\prime=\xi_{k,n}| {\mathcal S}\bigr)= \tfrac{1}{2}(1+e^{-2|u_2-u_1|})$.
Analogously we may define a stochastic flow of   kernels on ${\mathbf R}$ from a pair of $\theta$-coupled Brownian webs $\bigl({\mathcal W}, {\mathcal W}^\prime\bigr)$ via, 
\begin{equation}
K_{s,t}(x,A)= {\mathbf P} \bigl ( W^\prime_{x,s}(t) \in A| {\mathcal W} \bigr)
\end{equation}
for any $A$, a Borel subset of ${\mathbf R}$, and 
where  $W_{x,s}^\prime(t)$ denotes the position at time $t\geq s$ of  the almost surely unique path in the web ${\mathcal W}^\prime$ which starts from $(x,s)$. We will call this the erosion flow with parameter $\theta$. 

A powerful approach to describing a flow of kernels $K$ on a space $E$ is by means of its family of  $N$-point motions. For each integer  $N \geq 1$, the $N$th point 
motion 
of the flow is a Markov process on $E^N$. Formally it may be described by means of its semigroup which is given by
\begin{equation}
P^N_t(x,A)= {\mathbf E}\bigl[ K_{0,t}(x_1,A_1)K_{0,t}(x_2,A_2) \ldots K_{0,t}(x_N,A_N) \bigr],
\end{equation}  
for all $x=(x_1,x_2, \ldots, x_N) \in E^N$ and $A=A_1\times A_2\times \cdots \times A_N \in {\mathcal E}^N$.
Informally  it should be thought of as describing  the motion of $N$  infinitesimal particles sampled from the flow of mass, or if $K$ is 
interpreted as a random environment governing the evolution of a particle, then take $N$ such particles and let them evolve conditionally 
independently given $K$.
Notice that the family of $N$ point motions is consistent in that  any $M$ co-ordinates (regardless of order) taken from the 
$N$-dimensional process are distributed as the $M$-dimensional process in the family. The law of the  flow $K$ (in the sense of  finite-dimensional 
distributions)  is uniquely determined by the associated family of  $N$-point semigroups $\bigl(P^N_t; t\geq 0 \bigr)$ for $N \geq 1$.

Turning to the erosion flow, we seek to identify its $N$-point motions. It is clear from the discrete counterpart that the $N$-point motion should be a diffusion in ${\mathbf R}^N$, whose $N$ co-ordinates are 
all Brownian  motions with  the only interactions  between co-ordinates being local:  when
co-ordinates differ in value they evolve independently of each other.
Le Jan and Raimond, \cite{MR2052863}, see also \cite {MR2092875}, first gave an example of a flow  whose $N$-point motions behaved like this. In \cite{hw} we described more geneal flows of this type, whose $N$-point motions were characterized  by means of martingale problems which involve a family of parameters $ \bigl( \theta(k:l); k,l \geq 0 \bigr)$. 
In  a certain excursion theoretic sense, $\theta(k:l)$,  for $k,l \geq 1$, describes the rate at which 
a group of $n=k+l$ particles moving together splits into two groups one  consisting of  $k$ particles and the other of $l$ particles. 
We will prove that the $N$-point motions of the erosion flow, constructed from $\theta$-coupled webs, correspond to parameters given by
\begin{equation}
\theta(1:1)= \theta, \;\; \theta(1:k)=\theta(k:1)=\theta/2 \text{ for } k \geq 2, \text{ and }  \theta(k:l)=0  \text{ if both } k,l \geq 2.
\end{equation}
Thus, a cluster of $n$ particles looses particles only one at a time. Interpreting the flow of kernels as a flow of mass, particles are infinitesimal, 
and this means mass is lost from a cluster in a continuous way rather than by splitting. This is analogous to the phenomenon known as erosion in fragmentation theory, \cite{bertoin}.

This paper is organized as follows. In Section 2, we recall the construction and characterization of the Brownian web given by Fontes et al., \cite{FINR}. In Section 3, we introduce and characterize $\theta$-coupled Brownian webs. In Section 4, we describe a method of approximating $\theta$-coupled webs by  pairs of coupled  webs in which there is high frequency switching between the two regimes: one in which paths in the two copies evolve independently of each other, and the other in which they  coalesce. Using these approximations, we prove in section 5 a Markov property for $\theta$-coupled webs, which is a  form of the semigroup property for the  transition semigroup of the  stationary process $\bigl( {\mathcal W}(u); u \in {\mathbf R}\bigr)$ mentioned above. We then return to the approximation of $\theta$-coupled webs by switching and strengthen the mode of convergence taking place. Finally in Section 6 of the paper we are able to identify the $N$-point motions of the flow of kernels constructed from $\theta$-coupled webs. Section 7 contains some comments on possible generalisations.  

\section{Preliminaries on the Brownian web}

We begin with a precise description of the Brownian web.   Let 
$\bigl(\bar{{\mathbf R}}^2, \rho\bigr)$  be  the completion of ${\mathbf R}^2$ with respect to the metric $\rho$ which is given by 
\begin{equation}
\rho \bigl( (x_1,t_1), (x_2,t_2)\bigr)=\max\left\{\left| \frac{\tanh(x_1)}{1+|t_1|}-\frac{\tanh(x_2)}{1+|t_2|}\right|, | \tanh(t_1)-\tanh(t_2)| \right\}.
\end{equation}
This space plays the role of space-time. As a set we may identify $\bar{{\mathbf R}}^2$ with ${\mathbf R}^2 \cup \{(\pm \infty, t ): t \in {\mathbf R} \} \cup \{( *, +\infty), (*, -\infty)\}$.
 
Next  we construct a space of paths with specified starting times. For $t_0 \in [-\infty,\infty]$ let ${\mathcal C}[t_0]$ denote the set of functions $f$ from $[t_0,\infty]$ to $[-\infty, +\infty]$ satisfying  $f(+\infty)=(*,+\infty)$ and $f(-\infty)=(*,-\infty)$ in the case $t_0=-\infty$, such that $t \mapsto (t,f(t))$ is continuous with respect to the metric $\rho$. Then we define the space $\Pi$ to be
\begin{equation}
 \Pi= \bigcup_{t_0 \in [-\infty,+\infty]} {\mathcal C}[t_0] \times \{t_0\}.
\end{equation}
Endowed with the metric
\begin{multline}
d \bigl( (f_1,t_1), (f_2,t_2)\bigr)= \\
\max\left\{ \sup_{t \geq t_1 \wedge t_2} \left| \frac{\tanh(f_1( t \vee t_1))}{1+|t|}-\frac{\tanh(f_2(t \vee t_2))}{1+|t|}\right|, | \tanh(t_1)-\tanh(t_2)| \right\},
\end{multline}
$\Pi$ becomes a complete, separable, metric space. 

Finally we define ${\mathcal H}$ to be  the space of compact subsets of $\Pi$ and let $d_{\mathcal H}$ be the induced Hausdorff metric. 
Let ${\mathcal B}_{\mathcal H}$ denote the Borel $\sigma$-algebra on ${\mathcal H}$ associated with $d_{\mathcal H}$.

\begin{theorem}[Fontes et al., \cite{FINR}]
\label{thm:finr}
There exists a $\bigl({\mathcal H},{\mathcal B}_{\mathcal H} \bigr)$-valued random variable ${\mathcal W}$, called the Brownian web, whose distribution is uniquely determined by the following three properties.
\begin{enumerate}
\item From any deterministic point $(x,t) \in {\mathbf R}^2$ there is almost surely a unique path $W_{x,t}$ starting from $(x,t)$.
\item For any deterministic list of $n$ points, $(x_1,t_1), (x_2,t_2) \ldots, (x_n,t_n)$, the joint distribution of $W_{x_1,t_1}, W_{x_2,t_2}, \ldots,W_{x_n,t_n}$ is that of a coalescing system of Brownian motions.
\item For any  deterministic, countable, dense subset ${\mathcal D}$ of ${\mathbf R}^2$, almost surely, ${\mathcal W}$ is the closure in $( \Pi,d)$ of
$\{ W_{x,t}: (x,t) \in {\mathcal D} \}$.
\end{enumerate}
\end{theorem}

Recall a pair of processes  $\bigl(X_1(t),X_2(t);  t \geq 0 \bigr)$ is said to be a pair of coalescing Brownian motions if both $X_1,X_2$ are Brownian motions, relative to a common filtration, and their bracket satisfies
\begin{equation}
\langle X_1, X_2 \rangle (t)= ( t-T)^+,
\end{equation}
where $T= \inf\{ t \geq 0: X_1(t)=X_2(t)\}$. Informally paths evolve as independent Brownian motions until the instant they first meet, at which point they coalesce, and  thereafter they evolve  identically. This notion can easily be extended to a system of several coalescing paths, with different starting times. 

 Suppose that the web ${\mathcal W}$ is defined on an underlying complete probability space $\bigl(\Omega,{\mathcal F}, {\mathbf P} \bigr)$. All sub $\sigma$-algebras of  ${\mathcal F}$ are  assumed to contain all null events. For $-\infty  < t < +\infty$ let ${\mathcal F}_{t}$ be the   sub-$\sigma$-algebra of ${\mathcal F}$ generated by  random variables  of the form  $W_{x,s_1}(s_2) $  with $s_1 \leq s_2 \leq t$. Clearly $\bigl({\mathcal F}_{t}; t\in {\mathbf R} \bigr)$ forms a filtration: it is the natural filtration of the web ${\mathcal W}$. In particular it can be shown without great difficulty that for each $(x,t) \in {\mathbf R}^2$ the  path  
 $\bigl(W_{x,t}(u); u \geq t\bigr)$  is a Brownian motion relative to $\bigl({\mathcal F}_{u}; u \geq t  \bigr)$. More generally for $-\infty \leq s < t\leq \infty $ we define ${\mathcal F}_{s,t}$ to be the   sub-$\sigma$-algebra of ${\mathcal F}$ generated by   random variables of the form $W_{x,u_1}(u_2) $  with $s \leq u_1 \leq u_2 \leq t$.
 
\begin{proposition}
\label{prop:ctsprod}
For any $s<t<u$ we  have
\begin{itemize}
\item[(1)] ${\mathcal F}_{s,t}$ and ${\mathcal F}_{t,u}$ are independent;
\item[(2)] ${\mathcal F}_{s,t}$ and ${\mathcal F}_{t,u}$ together generate ${\mathcal F}_{s,u}$.
\end{itemize} 
\end{proposition}

\begin{proof} 
A system of coalescing Brownian motions possesses the Markov property. In particular suppose $\bigl(W_1, W_2, \ldots W_n)$ is such a system with $W_k$ starting from $(x_k,t_k)$, and $t_1, t_2, \ldots t_m <t $ and $t_{m+1}, t_{m+2}, \ldots t_n \geq t$ for some given $t \in {\mathbf R}$. Then $ \bigl(W_{m+1}, W_{m+2}, \ldots, W_{n} \bigr)$ is independent of the evolution of $\bigl(W_1, W_2, \ldots W_m \bigr)$ prior to time $t$. Property (1) follows from this.

Property (2) follows from  the flow property. Let us write $W(s,t,x)$ for the random variable $W_{x,s}(t)$. Then as in Lemma 8.3 of \cite{tw}, for any $x \in {\mathbf R}$, and $ s \leq t \leq u$, almost surely,
\[
W(s,u,x)= W(t,u,W(s,t,x)).
\]
Implicit here is the statement that almost surely there is a unique path starting from $W(s,t,x)$. 
\end{proof} 
 
By virtue  of the properties asserted by the preceding proposition $\bigl( \Omega,{\mathcal F}, {\mathbf P} \bigr)$ equipped with the $\sigma$-algebras $\bigl\{{\mathcal F}_{s,t}; s<t\bigr\}$ 
is called a continuous product of probability spaces, see Tsirelson \cite{tsirelson}. We will refer to $\bigl\{{\mathcal F}_{s,t}; s<t\bigr\}$  as the factorization generated by the web ${\mathcal W}$.

\section{Coupled webs}

We begin with repeating from the introduction the definition of $\theta$-coupled Brownian motions. A pair of Brownian motions $\bigl(X(t); t \geq 0\bigr)$ and $\bigl(X^\prime(t); t \geq 0\bigr)$ defined on a common 
probability space 
are  $\theta$-coupled, where $\theta$ is  a  positive real parameter, if, $X$ and $X^\prime$ are both Brownian motions relative to some 
common filtration, and,
\begin{align}
\label{quad}
\langle X,X^\prime\rangle (t)=\int_0^t {\mathbf 1}_{( X(s)=X^\prime(s))} ds  \qquad &  t\geq 0, \\
\label{locall}
L^0_t\bigl( X-X^\prime\bigr) = 2\theta \int_0^t {\mathbf 1}_{( X(s)=X^\prime(s))} ds  \qquad& t\geq 0.
\end{align}

\begin{proposition}
\label{bm:unique} 
For each given starting point $ (x_1,x_2) \in {\mathbf R}^2$, and  parameter $\theta>0$,  there exists a pair of $\theta$-coupled Brownian motions starting from $(x_1,x_2)$ and its law is uniquely determined.  
\end{proposition}
\begin{proof} This is a special case of Proposition \ref{prop:gencouple} a proof of which is given in Section 7.
\end{proof}

We may extend the notion of $\theta$-coupled Brownian motions to families of paths as follows.  An ${\mathbf R}^{n+m}$-valued process $ \bigl( X_1(t), X_2(t), \ldots X_m(t), X^\prime_1(t), X_2^\prime(t), \ldots, X_n^\prime(t); t \geq 0 \bigr)$ is a sticky-coalescing system with parameter $\theta$  if the following properties hold.
\begin{itemize}
\item[(1)] Each of the $(m+n)$ processes $X_1,X_2, \ldots X_m$ and $X^\prime_1,X^\prime_2, \ldots X^\prime_n$ is a  Brownian motion relative to some common filtration.
\item[(2)] For $i,j \in \{1,2, \ldots,m \} $,  $X_i$ and $X_j$ are a pair of coalescing Brownian motions.  Similarly, for $i,j \in \{1,2, \ldots,n \} $,  $X^\prime_i$ and $X^\prime_j$ are a pair of coalescing 
Brownian motions.
\item[(3)]
For $i\in \{1,2, \ldots,m \} $ and $j \in \{1,2, \ldots,n \} $, $X_i$ and $X^\prime_j$ are a pair of $\theta$-coupled Brownian motions.
\end{itemize}
Standard localization arguments, allow us to deduce from Proposition \ref{bm:unique}, that for  given any point $(x_1,x_2, \ldots x_m, x_1^\prime,x_2^\prime, \ldots, x_n^\prime)$  there exists a sticky-coalescing  system with parameter $\theta$ starting from this point, and it possesses a  uniquely determined law. 
By virtue of this uniqueness and their definition, it is clear that the   laws of sticky-coalescing systems (SCS) of Brownian motions have a natural consistency property: if an ${\mathbf R}^{d}$-valued process is a  SCS then so too is any ${\mathbf R}^{d^\prime}$-valued process, $ d^\prime\leq d$, formed by taking some subset of the  coordinates  of the given process.

We can generalize the definition of sticky-coalescing systems by allowing the paths to have different starting times. In particular a pair of processes $\bigl(X, X^\prime)$ with $X$ indexed by $ t \in [s_1,\infty)$ and and $X^\prime$ indexed $t \in [s_2,\infty)$ where $s_1, s_2 \in {\mathbf R}$ may be  arbitrary will be called a pair of $\theta$-coupled Brownian motions if there exists some filtration $\bigl({\mathcal F}_t; t \geq s_1 \wedge s_2 \bigr)$ such that  $\bigl(X_i(t); t \geq s_i)$ is a Brownian motion relative to $\bigl({\mathcal F}_t; t \geq s_i \bigr)$ for $i=1,2$, and such that
 the pair of processes $\bigl(X_1((s_1 \vee s_2)+t), X_2((s_1 \vee s_2)+t); t \geq 0 \bigr)$  satisfies \eqref{quad} and \eqref{locall}.
The definition of sticky-coalescing systems of $(m+n)$ paths generalizes in a similar fashion.

We turn now to the construction of coupled pairs of webs as  suggested in the Introduction.
An ${\mathcal H} \times {\mathcal H}$-valued random variable $\bigl({\mathcal W}, {\mathcal W}^\prime\bigr)$ is called a coupled pair of Brownian webs if 
both ${\mathcal W}$ and ${\mathcal W}^\prime$ are distributed  according to the law of the Brownian web.   Let $\bigl({\mathcal F}_t; t \in {\mathbf R} \bigr)$ and $\bigl({\mathcal F}^\prime_t; t \in {\mathbf R} \bigr)$ be the natural  filtrations of ${\mathcal W}$ and ${\mathcal W}^\prime$ respectively. 
The webs ${\mathcal W}$ and $ {\mathcal W}^\prime$ will be said to be co-adapted if
for each $(x,t) \in {\mathbf R}^2$ the paths  $\bigl(W_{x,t}(u); u \geq t \bigr)$ and $\bigl(W^\prime_{x,t}(u); u \geq t\bigr)$, starting from $(x,t)$, and  contained in ${\mathcal W}$ and ${\mathcal W}^\prime$ respectively,  are Brownian motions relative to the filtration $\bigl({\mathcal F}_u \vee {\mathcal F}_u^\prime; u \geq t)$.

\begin{theorem}
\label{coupledwebs}
There exists a $\bigl({\mathcal H} \times{\mathcal H}, {\mathcal B}_{\mathcal H} \otimes  {\mathcal B}_\mathcal{H} \bigr)$-valued random variable $\bigl({\mathcal W},{\mathcal W}^\prime\bigr)$ defined on some probability space $\bigl(\Omega, {\mathcal F}, {\mathbf P} \bigr)$ whose law is uniquely determined by the following properties.
\begin{itemize}
\item[(1)] ${\mathcal W}$ and ${\mathcal W}^\prime$ are both distributed as the Brownian web.
\item[(2)] ${\mathcal W}$ and ${\mathcal W}^\prime$ are co-adapted.
\item[(3)] For any pair of deterministic points $(x,t)$ and $(x^\prime,t^\prime)$ in ${\mathbf R}^2$, the paths $W_{x,t}$ and $W^\prime_{x^\prime,t^\prime}$
are a $\theta$-coupled pair of Brownian motions.
\end{itemize}
\end{theorem}

\begin{proof} From the consistency of  sticky-coalescing systems  of Brownian motions  and   Kolmogorov's extension theorem,
 given two countable dense subsets ${\mathcal D}=\{ (x_1,t_1),(x_2,t_2), \ldots \}$ and ${\mathcal D}^\prime=\{ (x^\prime_1,t^\prime_1),(x^\prime_2,t^\prime_2), \ldots \}$ of ${\mathbf R}^2$,  we can assert the existence of an infinite family of processes 
  $\bigl(X_1, X_2, \ldots, X^\prime_1, X^\prime_2, \ldots \bigr)$ defined on a common probability space such that for any $m,n \geq 1$, the ${\mathbf R}^{m+n}$ valued process $\bigl(X_1,  \ldots, X_m, X^\prime_1, \ldots, X^\prime_n \bigr)$ is a SCS starting from  
  $\bigl((x_1,t_1), \ldots, (x_m,t_m),(x^\prime_1,t^\prime_1),\ldots,(x^\prime_n,t_n^\prime)\bigr)$. Each path $X_k$, or  $X^\prime_l$, may be treated as  a point in the space $\Pi$, and  then $\bigl\{ X_k: k \geq 1\big\}$ and  $\bigl\{ X^\prime_l: l \geq 1\big\}$ are each random subsets of $\Pi$. We define  ${\mathcal W}$ to be the closure of the former  and ${\mathcal W}^\prime$ to be the closure of the latter. As in the proof of Theorem \ref{thm:finr} given in \cite{FINR}, ${\mathcal W}$  and ${\mathcal W}^\prime$ are each Brownian webs.

We need to show that properties (2) and (3) hold for the pair of webs ${\mathcal W}$  and ${\mathcal W}^\prime$ just constructed. We verify (3) first. 
We are given two points $(x,t)$ and $(x^\prime,t^\prime)$. If $(x,t) \in {\mathcal D}$ and $(x^\prime,t^\prime) \in {\mathcal D}^\prime$ then the paths $W_{x,t}$ and $W^\prime_{x^\prime,t^\prime}$ are $\theta$-coupled Brownian motions  by construction. In general we may choose sequences $(x_n,t_n) \rightarrow (x,t)$ and   $(x^\prime_n,t^\prime_n) \rightarrow (x^\prime,t^\prime)$ with the points $(x_n,t_n) \in {\mathcal D}$ and $(x_n^\prime, t_n^\prime)\in {\mathcal D}^\prime$. Then on the one hand, a coupling argument, see \cite{tw},  shows that $W_{x_n,t_n} \rightarrow  W_{x,t}$ in $\Pi$ almost surely, and similarly $W_{x^\prime_n,t^\prime_n}^\prime \rightarrow  W^\prime_{x^\prime,t^\prime}$. On the otherhand, the law of a pair of $\theta$-coupled Brownian motions depends continuously on the starting points, and hence the law of the pair $\bigl(W_{x,t},W^\prime_{x^\prime,t^\prime}\bigr)$ is that of a pair of  $\theta$-coupled Brownian motions.

We verify (2) by arguing as follows.
Let $\bigl({\mathcal G}^{n}_t; t \in {\mathbf R} \bigr)$ be the filtration generated by   $\bigl( X_k, X^\prime_k; k=1,2,\ldots, n \bigr)$.  By definition of a SCS, if $n \geq m$ then $\bigl(X_m(t); t \geq t_m\bigr)$ is a Brownian motion relative to $\bigl({\mathcal G}^n_t; t \geq t_m \bigr)$, and likewise for $\bigl(X^\prime_m(t); t \geq t_m\bigr)$.
Let $\bigl( {\mathcal G}^\infty_t ; t \in {\mathbf R} \bigr)$ be the smallest filtration containing  $\bigl( {\mathcal G}^n_t ; t \in {\mathbf R} \bigr)$ for every $n$. Since  the path $W_{x,t}$ for arbitrary $(x,t)\in {\mathbf R}^2$ is a limit of paths $W_{x_n,t_n}$ with $(x_n,t_n) \in {\mathcal D}$,  it follows that $\bigl(W_{x,t}(u); u \geq t\bigr)$ is a Brownian motion relative to $\bigl( {\mathcal G}^\infty_u ;u \geq t \bigr)$  and similarly for $W^\prime_{x^\prime,t^\prime}$. But it must be that ${\mathcal G}^\infty_t= {\mathcal F}_t \vee {\mathcal F}_t^\prime$ where $\bigl({\mathcal F}_t; t \in {\mathbf R} \bigr)$ and  $\bigl( {\mathcal F}_t^\prime; t \in {\mathbf R} \bigr)$ are the natural filtrations of the webs ${\mathcal W}$ and ${\mathcal W}^\prime$, and so the webs are co-adapted.

To prove that the  uniqueness assertion, we observe that  if $\bigl({\mathcal W},{\mathcal W}^\prime\bigr)$ is any ${\mathcal H} \times {\mathcal H}$- valued random variable with the three given properties then the  joint distribution of the paths $ \bigl( W_{x_1,t_1},W_{x_2,t_2}, \ldots W_{x_m,t_m},
W^\prime_{x^\prime_1,t^\prime_1},W^\prime_{x^\prime_2,t^\prime_2}, \ldots W^\prime_{x^\prime_n,t^\prime_n}\bigr)$ is that of a SCS. Letting ${\mathcal D}$ and ${\mathcal D}^\prime$ be as before, we know from the property (3) of Theorem \ref{thm:finr} that ${\mathcal W}$ is the closure of the subset $\bigl\{ W_{x,t}; (x,t) \in {\mathcal D}\bigr\}$, and similarly for ${\mathcal W}^\prime$ for which it follows the distribution of $\bigl({\mathcal W},{\mathcal W}^\prime \bigr)$ is the same as that of the pair of webs constructed in the existence argument.
\end{proof}

The following property of $\theta$-coupled webs strengthens the notion of being co-adapted: in the terminology of Tsirelson it says that $\theta$-coupled webs are a joining of continuous products.

\begin{proposition}
\label{joining}
Let $\bigl\{{\mathcal F}_{s,t}; s<t \bigr\}$ and  $\bigl\{{\mathcal F}^\prime_{s,t}; s<t \bigr\}$ be the factorizations generated by Brownian webs ${\mathcal W}$ and ${\mathcal W}^\prime$ which are $\theta$-coupled. Then for $s<t<u$,
\[
{\mathcal F}_{s,t} \vee {\mathcal F}_{s,t}^\prime \text{ is independent of } {\mathcal F}_{t,u} \vee {\mathcal F}_{t,u}^\prime.
\]
\end{proposition}

\begin{proof}
Consider $t_1,  t_2 \ldots t_m, t_1^\prime, t_2^\prime, \ldots, t_n^\prime \geq t$ and $x_1, x_2 ,\ldots,x_,  \in {\mathbf R}$. Then since  webs ${\mathcal W}$ and ${\mathcal W}^\prime$ are co-adapted each path $\bigl(W_{x_k,t_k}(v); v \geq t_k \bigr)$  is a Brownian motion relative to the filtration $\bigl( {\mathcal F}_v \vee {\mathcal F}_v^\prime; v \geq t_k \bigr)$, and similarly for the path $W^\prime_{x^\prime_k,t^\prime_k}$. We may deduce from this and the characterization of SCS  that $ \bigl( W_{x_1,t_1},W_{x_2,t_2}, \ldots W_{x_m,t_m},
W^\prime_{x^\prime_1,t^\prime_1},W^\prime_{x^\prime_2,t^\prime_2}, \ldots W^\prime_{x^\prime_n,t^\prime_n}\bigr)$ is independent of ${\mathcal F}_t \vee {\mathcal F}_t^\prime$. The result follows, since ${\mathcal F}_{s,t} \vee {\mathcal F}_{s,t}^\prime$ is contained in ${\mathcal F}_t \vee {\mathcal F}_t^\prime$, and ${\mathcal F}_{t,u} \vee {\mathcal F}_{t,u}^\prime$ is generated by random variables of the form $W_{x_k,t_k}(u_k)$ and 
$W^\prime_{x^\prime_k,t^\prime_k}(u_k^\prime)$ with $t\leq t_k \leq u_k \leq  u$  and $t\leq t^\prime_k \leq u^\prime_k \leq  u$.
\end{proof}

\section{A convergence result}

We do not have sufficient tools to study $\theta$-coupled webs  directly from the characterization given in the previous section. Instead we must
use some approximation regime. It would be natural to use the discrete approximations suggested by the discussion in the introduction. But, in fact,
we use a different method, suggested by Tsirelson's theory of continuous products, in which the approximations are also coupled Brownian webs, which is technically advantageous, switching at high frequency between evolving independently of one another and evolving identically.   

We begin with a convergence result for  coupled Brownian paths.
Given a parameter $p\in[0,1]$ and an  integer $n \geq 1$ we define a $(p,n)$-coupling of Brownian motions as follows. Let  $(Y_k; k \geq 0 \bigr)$ be a sequence of independent Bernoulli($p$) random variables. Then conditionally on $(Y_k; k \geq 0 \bigr)$, the process $\bigl(X(t),X^\prime(t);t\geq 0\bigr)$ is a time-inhomogeneous diffusion, evolving as:
\begin{itemize}
\item[(i)]
a pair of independent Brownian motions whilst $t \in [k/n,(k+1)/n]$ if $Y_k=1$;
\item[(ii)]
a pair of coalescing Brownian motions whilst $t \in [k/n,(k+1)/n]$ if $Y_k=0$.
\end{itemize}

In the following convergence in distribution means weak convergence of probability measures on the path space $ {\mathbf C} \bigl( 
[0,\infty),{\mathbf R}^2 \bigr)$.

\begin{proposition}
\label{bm:con1}
Let $\bigl(X,X^{(n)}\bigr)$ be a sequence of $(p,n)$-coupled Brownian motions with starting points  not depending on $n$. Suppose that $p=p(n)$ 
satisfies
\[
\lim_{n \rightarrow \infty} \sqrt{\frac{n}{\pi}}{p(n)} = \theta \in (0,\infty). 
\]
Then, as $n$ tends to infinity, $\bigl(X,X^{(n)}\bigr)$ converges in distribution to a pair of $\theta$-coupled Brownian motions.
\end{proposition}
\begin{proof}
The sequence of laws of $\bigl(X,X^{(n)}\bigr)$ is tight since the  distribution of $X$ and $X^{(n)}$ evidently do not depend on $n$. Thus it is sufficient to show any subsequence $\bigl(X,X^{(n_k)}\bigr)$ which converges in distribution must converge to a pair of $\theta$-coupled Brownian motions. Assume for notational simplicity  that  $\bigl(X,X^{(n)}\bigr)$ itself converges in distribution to a process  $\bigl(X,X^\prime\bigr)$. 
It is clear that  $X$ and $X^\prime$ must each distributed as Brownian motions, and moreover $X$ and $X^\prime$ are each martingales relative to the natural filtration generated by $\bigl(X,X^\prime\bigr)$, since the corresponding statement holds for $X$ and $X^{(n)}$. So in order to appeal to the uniqueness assertion of Proposition \ref{bm:unique}, we must verify \eqref{quad} and \eqref{locall} hold.

Our first task is to determine the quadratic covariation of $X$ and $X^\prime$.
We begin with the observation that it is a general fact, valid for any continuous semimartingale $Z$, that
$ \int_0^t {\mathbf 1}(Z(s)=0) d \langle Z \rangle(s)=0$. Applying this to $Z=X-X^\prime$, and using $\langle X \rangle (t)= \langle X^\prime \rangle (t)=t$ we deduce that
\begin{equation}
\label{gen}
\int_0^t {\mathbf 1}(X(s)=X^\prime(s))d\langle X, X^\prime \rangle (s)=\int_0^t {\mathbf 1}(X(s)=X^\prime(s))ds.
\end{equation}
Consider two times $0\leq t_1 \leq t_2$. Let $g: {\mathbf C} \bigl( [0,\infty),{\mathbf R}^2 \bigr) \rightarrow {\mathbf R}$ be  
non-negative, bounded, continuous, and measurable with respect to ${\mathcal B}_{t_1}$, where $\bigl({\mathcal B}_t; t \geq 0 \bigr)$ is the 
filtration generated by the co-ordinate process. The mapping  $\alpha \mapsto \int_{t_1}^{t_2} {\mathbf 1}(\alpha_1(s)=\alpha_2(s))ds$  
is 
upper semicontinuous relative to the local uniform topology on $ {\mathbf C} \bigl( [0,\infty),{\mathbf R}^2 \bigr)$.  Thus by weak 
convergence,
\begin{multline}
\label{limsup}
{\mathbf E} \left[ g\bigl(X,X^\prime\bigr)\int_{t_1}^{t_2} {\mathbf 1}(X(s)=X^\prime(s))ds \right] \geq \\
\limsup_{n \rightarrow \infty} \; {\mathbf E} \left[ 
g\bigl(X,X^{(n)}\bigr)\int_{t_1}^{t_2} {\mathbf 1}(X(s)=X^{(n)}(s))ds \right].
\end{multline}Also, since  $ \langle X, X^{(n)} \rangle(t) = \int_{0}^{t} {\mathbf 1}(X(s)=X^{(n)}(s))ds $ for every $n$, we have,
\begin{multline}
\label{limsup2}
\lim_{n \rightarrow \infty} \; {\mathbf E} \left[ 
g\bigl(X,X^{(n)}\bigr)\int_{t_1}^{t_2} {\mathbf 1}(X(s)=X^{(n)}(s))ds \right]= \\
\lim_{n \rightarrow \infty} \; {\mathbf E} \left[ 
g\bigl(X,X^{(n)}\bigr)\bigl\{X(t_2)X^{(n)}(t_2)-X(t_1)X^{(n)}(t_1)\bigr\} \right] =
\\
{\mathbf E} \left[ 
g\bigl(X,X^\prime\bigr)\bigl\{X(t_2)X^\prime(t_2)-X(t_1)X^\prime(t_1)\bigr\} \right] ,
\end{multline}
the last equality holding by weak convergence and uniform integrability. Combining \eqref{limsup} and \eqref{limsup2}, since $g$, $t_1$ and $t_2$ are arbitrary, we deduce that
$ \langle X, X^\prime \rangle(t) - \int_{0}^{t} {\mathbf 1}(X(s)=X^\prime(s))ds $ must be a non-increasing process.  On the other hand, since \eqref{limsup2} also gives  
\begin{equation}
{\mathbf E} \left[ 
g\bigl(X,X^\prime\bigr)\bigl\{X(t_2)X^\prime(t_2)-X(t_1)X^\prime(t_1)\bigr\} \right] \geq 0,
\end{equation}
we must have $ \langle X, X^\prime \rangle(t)$ is a non-decreasing process. In view of  \eqref{gen} we thus may deduce that $ \langle X, X^{\prime} \rangle(t) = \int_{0}^{t} {\mathbf 1}(X(s)=X^\prime(s))ds $ as desired.
Moreover, it then follows from \eqref{limsup2} that
\begin{multline}
\label{lim}
\lim_{n \rightarrow \infty} \; {\mathbf E} \left[ 
g\bigl(X,X^{(n)}\bigr)\int_{t_1}^{t_2} {\mathbf 1}(X(s)=X^{(n)}(s))ds \right]=
\\
{\mathbf E} \left[ 
g\bigl(X,X^\prime\bigr)\int_{t_1}^{t_2} {\mathbf 1}(X(s)=X^\prime(s))ds \right],
\end{multline}
which will be useful to us shortly.

We turn now to computing the local time at zero of $X- X^\prime$. For $t\geq 0$ define $\kappa_t(x)={\mathbf E} \bigl[ |B(2t)+x| \bigr]-|x|$, where $B$ is a standard Brownian motion starting from zero.
By Brownian scaling we have  $\kappa_t(x)= t^{1/2}\kappa_1(t^{-1/2}x)$. Using this, elementary calculations show that
\begin{equation}
\label{kappalimit}
t^{-1/2} \kappa_t(x) \downarrow \frac{2}{\sqrt{\pi}}{\mathbf 1}(x=0) \quad \text{ as }  t \downarrow 0.
\end{equation}
 Let $\bigl({\mathcal F}^{(n)}_{t}; t \geq 0 \bigr)$ denote the natural filtration of the coupled Brownian motions $\bigl(X,X^{(n)}\bigr)$.
We observe that from the definition of  $\bigl(X,X^{(n)}\bigr)$ it follows easily that 
\begin{equation}
\label{dismart1}
| X(k/n)- X^{(n)}(k/n)| - p(n) \sum_{r=0}^{k-1} \kappa_{1/n} \bigl(X(r/n)- X^{(n)}(r/n)\bigr)
\end{equation}
for $k=0,1,2,\ldots$, defines a discrete parameter martingale relative to the filtration $\bigl({\mathcal F}^{(n)}_{k/n}; k\geq 0 \bigr)$.
As we noted above, $t^{-1/2} \kappa_t(x)$ decreases as $t$ decreases with $x$ fixed, and consequently, if $m \leq n$, then 
\[
| X(k/n)- X^{(n)}(k/n)| - p(n) \sqrt{\frac{m}{n}} \sum_{r=0}^{k-1} \kappa_{1/m} \bigl( X(r/n)- X^{(n)}(r/n)\bigr)
\]
defines a supermartingale. Let $g$, $t_1$ and $t_2$ be as before, and choose $t_1^n$ and $t_2^n$ so that $nt_1^n$ and $nt_2^n$ are integers with
$t_1^n \downarrow t_1$, and $t_2^n \downarrow t_2$ as $n$ tends to infinity. Then, for $m \leq n$,
\begin{multline*}
{\mathbf E} \bigl[ g(X,X^{(n)})\bigl\{ | X(t_2^n)- X^{(n)}(t_2^n)|-| X(t_1^n)- X^{(n)}(t_1^n)| \bigr\}\bigr] \leq 
\\
p(n) \sqrt{\frac{m}{n}} {\mathbf E} \left [  g(X,X^{(n)}) \left\{ \sum_{r=nt_1^n}^{nt_2^n-1} \kappa_{1/m} \bigl( X(r/n)- X^{(n)}(r/n)\bigr) \right\} \right]. 
\end{multline*}
By weak convergence the lefthandside converges to 
\[
{\mathbf E} \bigl[ g(X,X^{\prime})\bigl\{ | X(t_2)- X^{\prime}(t_2)|-| X(t_1)- X^{\prime}(t_1)| \bigr\}\bigr],
\]
and the righthandside to
\[
\theta \sqrt{m \pi}{\mathbf E} \left [  g(X,X^{\prime})  \int_{t_1}^{t_2} \kappa_{1/m} \bigl( X(s)- X^{\prime}(s)\bigr)ds \right].
\]
Thus
\[
| X(t)- X^\prime(t)| -  \theta \sqrt{m \pi} \int_0^t \kappa_{1/m} \bigl( X(s) -X^\prime(s) \bigr)ds,
\]
is a supermartingale relative to the filtration generated by $X$ and $X^\prime$. Letting $m$ tend to infinity, and appealing to \eqref{kappalimit}, we deduce that
\begin{equation}
\label{sub}
| X(t)- X^\prime(t)| - 2 \theta \int_0^t {\mathbf 1} \bigl( X(s) = X^\prime(s) \bigr)ds
\end{equation}
is a supermartingale too. To complete the proof we must show the same process is also a submartingale, then \eqref{locall} will follow by Tanaka's formula.

Let $\lambda_t(x)= {\mathbf E} \left[ \int_0^t {\mathbf  1} (Z(s)=0)ds \right]$, where $Z$ is distributed as the difference of two coalescing Brownian motions started at a distance  $x$ apart. We may represent $\lambda_t(x)$ as
\[
\int_0^t {\mathbf P}\bigl( T_x \in ds\bigr) (t-s),
\]
where $T_x$ is distributed as the first time the two Brownian motions meet. The quantity $\kappa_t(x)$ has a similar representation as
\[
\int_0^t {\mathbf P}\bigl( T_x \in ds\bigr)\sqrt{\frac{4(t-s)}{\pi }},
\]
and, comparing the two, we see that
\begin{equation}
\label{comp}
\sqrt{\frac{4}{\pi t}} \lambda_t(x) \leq \kappa_t(x).
\end{equation}
Now as a consequence of the construction of $(p,n)$-coupling we have,
\[
\int_0^{k/n} {\mathbf  1} (X(s)=X^{(n)}(s))ds - (1-p(n)) \sum_{r=0}^{k-1} \lambda_{1/n} \bigl(X(r/n)- X^{(n)}(r/n)\bigr)
\]
for $k=0,1,2,\ldots$, defines a discrete parameter martingale relative to the filtration $\bigl({\mathcal F}^{(n)}_{k/n}; k\geq 0 \bigr)$.
Subtracting a suitable multiple of   this from  the martingale given \eqref{dismart1},  and using \eqref{comp}, we deduce that
\[
| X(k/n)- X^{(n)}(k/n)|- \frac{2p(n)}{1-p(n)}\sqrt{\frac{n}{\pi}}\int_0^{k/n} {\mathbf  1} (X(s)=X^{(n)}(s))ds 
\]
is a submartingale.
From this, by a weak convergence argument as above together with \eqref{lim}, we may deduce that the process at \eqref{sub} is a submartingale as was required.
\end{proof}

The following extension  of the notion of a $(p,n)$-coupling will pay a pivotal role in proving the Markovian property of $\theta$-couplings in the 
next section.
Given parameters $p\in[0,1]$,  $\theta\in (0, \infty)$, and  integer  $n \geq 1$ we construct what we shall call a  $(p,\theta,n)$-coupling of Brownian motions as follows.
Let  $(Y_k; k \geq 0 \bigr)$ be a sequence of independent Bernoulli($p$) random variables. Then conditionally on $(Y_k; k \geq 0 \bigr)$, the process $\bigl(X_t,X^\prime_t; t \geq 0 \bigr)$ is a time-inhomogeneous diffusion, evolving as 
\begin{itemize}
\item[(i)]
a pair of independent Brownian motions whilst $t \in [k/n,(k+1)/n]$ if $Y_k=1$;
\item[(ii)]
a pair of $\theta$-coupled Brownian motions whilst $t \in [k/n,(k+1)/n]$ if $Y_k=0$.
\end{itemize}

\begin{proposition}
\label{bm:con2}
Let $\bigl(X,X^{(n)}\bigr)$ be a sequence of $(p,\theta_1,n)$-coupled Brownian motions with starting points  not depending on $n$. Suppose that
$p=p(n)$ satisfies
\[
\lim_{n \rightarrow \infty} \sqrt{\frac{n}{\pi}}{p(n)} = \theta_2 \in (0,\infty). 
\]
Then, as $n$ tends to infinity, $\bigl(X,X^{(n)}\bigr)$ converges in distribution to a pair of $(\theta_1+\theta_2)$-coupled Brownian motions.
\end{proposition}

\begin{proof}
The proof of the previous proposition applies verbatim, until we come to calculation of the local time of $X-X^\prime$ which requires some minor changes. 

The martingale at \eqref{dismart1} must be replaced with
\begin{multline}
| X(k/n)- X^{(n)}(k/n)| - 2\theta_1  \int_0^{k/n}{\mathbf 1}( X(s)=X^{(n)}(s)) ds \\
 - p(n) \sum_{r=0}^{k-1} \kappa_{1/n} \bigl(X(r/n)- X^{(n)}(r/n)\bigr).
\end{multline}
From this we deduce, using weak convergence as before,  together with \eqref{lim}, that 
\begin{equation}
\label{martsum}
| X(t)- X^{\prime}(t)| - 2(\theta_1+\theta_2) \int_0^{t}{\mathbf 1}( X(s)=X^{\prime}(s)) ds
\end{equation}
is a supermartingale. 

Let $\lambda^\theta_t(x)= {\mathbf E} \left[ \int_0^t {\mathbf 1} (Z(s)=0)ds \right]$, where $Z$ is distributed as the difference of two $\theta$-coupled Brownian motions started at a distance  $x$ apart. 
Then we deduce from the construction of $(p, \theta, n)$-coupling that,
\[
\int_0^{k/n} {\mathbf  1} (X(s)=X^{(n)}(s))ds - (1-p(n)) \sum_{r=0}^{k-1} \lambda^\theta_{1/n} \bigl(X(r/n)- X^{(n)}(r/n)\bigr)
\]
for $k=0,1,2,\ldots$, defines a discrete parameter martingale. Clearly $\lambda^\theta_t(x) \leq \lambda_t(x)$, and  comparing  the latter with $\kappa_t(x)$ as before we deduce that
\[
| X(k/n)- X^{(n)}(k/n)| - \left\{ 2\theta_1 +\frac{2p(n)}{1-p(n)}\sqrt{\frac{n}{\pi}} \right\} \int_0^{k/n}{\mathbf 1}( X(s)=X^{(n)}(s)) ds,
\]
is a submartingale. Hence by a now familiar  convergence argument, \eqref{martsum} defines a martingale, and from Tanaka  $L^0_t(X-X^\prime)= 2(\theta_1+\theta_2) \int_0^{t}{\mathbf 1}( X(s)=X^\prime(s)) ds$ as required.

\end{proof}

As we have seen before, in studying the web, pairs of  processes indexed by $ u \in [t_1,\infty)$ and $u \in [t_2,\infty)$ arise, where $t_1, t_2 \in {\mathbf R}$ may be  arbitrary. We may trivially extend the definition of  $(p,n)$ and  $(p,\theta,n)$-coupled  Brownian motions to cover such processes, meaning  that their joint evolution on each interval of the form $[\max(t_1,t_2) ,\infty) \cap [k/n, (k+1)/n]$ is given by a randomly chosen regime. The two previous propositions then still hold for such processes.   

Mimicking the construction of Theorem \ref{coupledwebs} we may construct coupled systems of coalescing paths and hence  obtain  $(p,n)$-coupled  and $(p,\theta,n)$-coupled webs $\bigl( {\mathcal W},{\mathcal W}^{\prime}\bigr)$ such that  for any pair of deterministic points $(x,t)$ and $(x^\prime,t^\prime)$ in ${\mathbf R}^2$, the paths $W_{x,t}$ and $W^{\prime}_{x^\prime,t^\prime}$
are a $(p,n)$-coupled (respectively $(p,\theta, n)$-coupled ) pair of Brownian motions.  Letting $\bigl\{{\mathcal F}_{s,t} \bigr\}$ and $\bigl\{{\mathcal F}^\prime_{s,t} \bigr\}$ denote the factorizations generated by ${\mathcal W}$ and ${\mathcal W}^\prime$,
the joint distribution of such  coupled webs is characterized by the fact that for any two sequences of   random variables $\bigl( \Phi_k( {\mathcal W}); k \in {\mathbf Z}\bigr)$ and $\bigl( \Psi_k( {\mathcal W}^\prime); k \in {\mathbf Z}\bigr)$  having values in $[0,1]$ and with $\Phi_k({\mathcal W})$ measurable with respect to ${\mathcal F}_{k/n,(k+1)/n}$ and  $\Psi_k({\mathcal W}^\prime)$ measurable with respect to ${\mathcal F}^\prime_{k/n,(k+1)/n}$,
\begin{equation}
\label{char}
{\mathbf E} \left[ \prod \Phi_k( {\mathcal W} )\Psi_k ({\mathcal W}^{\prime}) \right] = \prod_k \Bigl\{ (1-p){\mathbf E} \left[ \Phi_k( {\mathcal W}^1) \Psi_k ({\mathcal W}^2) \right]+ p{\mathbf E} \left [ \Phi_k( {\mathcal W}^3) \Psi_k ({\mathcal W}^4) \right] \Bigr\}.
\end{equation}
where on the righthandside $\bigl({\mathcal W}^3, {\mathcal W}^4\bigr)$ denotes a pair of independent  Brownian webs, and  $\bigl({\mathcal W}^1, {\mathcal W}^2\bigr)$ denotes a pair of  webs, either satisfying ${\mathcal W}^1={\mathcal W}^2$ with probability one, or in the case of a $(p,\theta,n)$-coupling, a pair of $\theta$-coupled webs. It is not difficult to see that $(p,n)$- and $(p,\theta,n)$-coupled 
webs are co-adapted.

\begin{theorem}
\label{web:con}
Let $\bigl({\mathcal W},{\mathcal W}^{(n)}\bigr)$ be a sequence of $(p,n)$-coupled Brownian webs. Suppose that $p=p(n)$ satisfies
\[
\lim_{n \rightarrow \infty} \sqrt{\frac{n}{\pi}}{p(n)} = \theta \in (0,\infty). 
\]
Then, as $n$ tends to infinity, $\bigl({\mathcal W},{\mathcal W}^{(n)}\bigr)$ converges in distribution to a pair of $\theta$-coupled Brownian webs.

More generally let $\bigl({\mathcal W},{\mathcal W}^{(n)}\bigr)$ be a sequence of $(p,\theta_1,n)$-coupled Brownian webs. Suppose that
\[
\lim_{n \rightarrow \infty} \sqrt{\frac{n}{\pi}}{p(n)} = \theta_2 \in (0,\infty). 
\]
Then, as $n$ tends to infinity, $\bigl({\mathcal W},{\mathcal W}^{(n)}\bigr)$ converges in distribution to a pair of $(\theta_1+\theta_2)$-coupled Brownian webs.
\end{theorem}
\begin{proof}
The laws of the pairs  $\bigl({\mathcal W},{\mathcal W}^{(n)}\bigr)$ are relatively compact, so it is sufficient to verify that any limit point is the
law of a pair of $\theta$-coupled webs ( respectively  a pair of $(\theta_1+\theta_2)$-coupled Brownian webs).  Suppose $\bigl({\mathcal W},{\mathcal W}^{(n)}\bigr)$ converges in distribution to $\bigl({\mathcal W}, {\mathcal W}^\prime\bigr)$.
Recall the result sometimes called Slutsky's lemma, see \cite{tsirelson}, that applies because the marginals  distributions of $\bigl({\mathcal W}, {\mathcal W}^\prime\bigr)$ do not depend on $n$:  
\[
\lim_{n \rightarrow\infty}
{\mathbf E} \bigl[ \Phi({\mathcal W}) \Psi({\mathcal W}^{(n)})\bigr]=
{\mathbf E} \bigl[ \Phi({\mathcal W}) \Psi({\mathcal W}^\prime ) \bigr],
\]
for any bounded measurable functions $\Phi$ and $\Psi$ defined on ${\mathcal H}$, there being no continuity requirement. 

  We must check that properties (1), (2) and (3) of Theorem \ref{coupledwebs} hold. (1) holds trivially.  For (2)  we must show that, for any $(x,t)\in {\mathbf R}^2$ the paths ${\mathcal W}_{x,t}$ and ${\mathcal W}_{x,t}^\prime$, are martingales with respect to the filtration $\bigl\{{\mathcal F}_u \vee {\mathcal F}^\prime_u; u \geq t \bigr\}$ generated by the two webs.  Consider  times $v>u>t$, and bounded measurable functions $\Phi$ and $\Psi$, each defined on ${\mathcal H}$, and with $\phi({\mathcal W})$ being ${\mathcal F}_t$-measurable and $\psi({\mathcal W}^\prime)$ being ${\mathcal F}^\prime_t$-measurable
Then, since ${\mathcal W}$ and ${\mathcal W}^{(n)}$ are co-adapted, for any real $\alpha$, 
\begin{multline*}
{\mathbf E} \bigl[ \Phi({\mathcal W}) \Psi({\mathcal W}^\prime )\exp \{i \alpha ({W}_{x,t}(v)-{ W}_{x,t}(u)\} \bigr]=
\\ \lim_{n \rightarrow\infty}
{\mathbf E} \bigl[ \Phi({\mathcal W}) \Psi({\mathcal W}^{(n)}) \exp \{i \alpha ({W}_{x,t}(v)-{W}_{x,t}(u)\} \bigr]=
\\
\lim_{n \rightarrow\infty}
{\mathbf E} \bigl[ \Phi({\mathcal W}) \Psi({\mathcal W}^{(n)})\bigr] \exp\{ -\alpha(v-u)^2/2\}
={\mathbf E}\bigl[ \Phi({\mathcal W}) \Psi({\mathcal W}^{\prime})\bigr] \exp\{ -\alpha(v-u)^2/2\}.
\end{multline*}
This shows that ${ W}_{x,t}$ is a  Brownian motion with respect to $\bigl\{{\mathcal F}_u \vee {\mathcal F}^\prime_u; u \geq t \bigr\}$, and a similar argument applies to ${ W}_{x,t}^\prime$ also.

Turning to property (3), consider any pair of deterministic points $(x,t)$ and $(x^\prime,t^\prime)$ in ${\mathbf R}^2$. Then, appealing to Slutsky's lemma again,   the pairs of  paths $\bigl(W_{x,t},W^{(n)}_{x^\prime,t^\prime}\bigr)$ converge in distribution to $\bigl(W_{x,t},W^{\prime}_{x^\prime,t^\prime}\bigr)$. That the latter is a pair of $\theta$-coupled Brownian motions then follows from the previous two propositions, or rather their generalization to the case of paths with unequal starting times.
\end{proof}

\section{The Markov property}

\begin{theorem}
\label{markov}
Suppose that ${\mathcal W}$,${\mathcal W}^\prime$ and ${\mathcal W}^{\prime\prime}$ are three copies of the Brownian web defined on a common probability space and satisfying
\begin{itemize}
\item[(i)]  $\bigl({\mathcal W},{\mathcal W}^\prime\bigr)$ is a $\theta_1$-coupling, and $\bigl({\mathcal W},{\mathcal W}^{\prime\prime}\bigr)$ is a $\theta_2$-coupling;
\item[(ii)] ${\mathcal W}^\prime$ and ${\mathcal W}^{\prime\prime}$ are conditionally independent given ${\mathcal W}$.
\end{itemize}
Then $\bigl({\mathcal W}^\prime,{\mathcal W}^{\prime\prime}\bigr)$ is a $(\theta_1+\theta_2)$-coupling of Brownian webs.
\end{theorem}

\begin{proof}
Let $p(n)$ vary with $n$ so that $\lim_{n \rightarrow \infty} \sqrt{\frac{n}{\pi}}{p(n)} = \theta_2$.
For $n \geq 1$ consider a triple of Brownian webs $\bigl({\mathcal W},{\mathcal W}^\prime,{\mathcal W}^{(n)}\bigr)$ defined on a common probability  space with 
$\bigl({\mathcal W},{\mathcal W}^\prime\bigr)$ being a $\theta_1$-coupling, $\bigl({\mathcal W},{\mathcal W}^{(n)}\bigr)$ being  a $(p(n),n)$-coupling, and ${\mathcal W}^\prime$ and ${\mathcal W}^{(n)}$ being conditionally independent given ${\mathcal W}$.
We observe using the characterization at \eqref{char}  that $\bigl( {\mathcal W}^\prime, {\mathcal W}^{(n)}\bigr)$ is a $(p(n),\theta_1,n)$-coupling.

The triple $\bigl({\mathcal W},{\mathcal W}^\prime,{\mathcal W}^{(n)}\bigr)$ converges in distribution to the triple $\bigl({\mathcal W},{\mathcal W}^\prime,{\mathcal W}^{\prime\prime}\bigr)$ defined in the statement of the theorem, for by Theorem \ref{web:con} and Slutsky's lemma,
\begin{multline*}
{\mathbf E} \bigl[ \Phi({\mathcal W}) \Psi ({\mathcal W}^\prime) \Upsilon({\mathcal W}^{(n)})\bigr]= 
{\mathbf E} \bigl[ \Xi({\mathcal W})  \Upsilon({\mathcal W}^{(n)})\bigr] 
\rightarrow {\mathbf E} \bigl[ \Xi({\mathcal W})  \Upsilon({\mathcal W}^{\prime\prime})\bigr]\\
={\mathbf E} \bigl[ \Phi({\mathcal W}) \Psi ({\mathcal W}^\prime) \Upsilon({\mathcal W}^{\prime\prime})\bigr],
\end{multline*}
where $ \Xi({\mathcal W})= {\mathbf E} \bigl[\Psi ({\mathcal W}^\prime)|{\mathcal W} \bigr]  \Phi({\mathcal W}) $.
 But with  a second appeal to Theorem \ref{web:con} we  also deduce that $\bigl({\mathcal W}^\prime,{\mathcal W}^{\prime\prime}\bigr)$ is a $(\theta_1+\theta_2)$-coupling, as required.
\end{proof}

With this Markov property available to us, we revisit the convergence of $(p,n)$-coupled webs to $\theta$-coupled webs.

\begin{proposition}
\label{prop:strong}
Let $\bigl({\mathcal W},{\mathcal W}^{(n)}\bigr)$ be a sequence of $(p,n)$-coupled Brownian webs. Suppose that $p=p(n)$ satisfies
\[
\lim_{n \rightarrow \infty} \sqrt{\frac{n}{\pi}}{p(n)} = \theta \in (0,\infty). 
\]
Let $ \bigl({\mathcal W}, {\mathcal W}^\prime\bigr)$ be a pair of $\theta$-coupled webs.
Then,  as $n$ tends to infinity, for any $\Phi$ defined on ${\mathcal H}$ with ${\mathbf E} \bigl[ \Phi({\mathcal W})^2\bigr]<\infty$,
\[
{\mathbf E}\bigl[ \Phi({\mathcal W}^{(n)}) | {\mathcal W} \bigr]  \rightarrow {\mathbf E}\bigl[ \Phi({\mathcal W}^\prime) | {\mathcal W} \bigr]  \text{ in } L^2.
\]
\end{proposition}
\begin{proof}
For $n \geq 1$, 
let $\bigl({\mathcal W}, {\mathcal W}^{(n)}, \tilde{{\mathcal W}}^{(n)} \bigr)$ be a triple of Brownian webs defined on a common probability space 
such that ${\mathcal W}^{(n)}$ and $\tilde{{\mathcal W}}^{(n)}$ are conditionally independent given ${\mathcal W}$ and such that both $\bigl({\mathcal W}, {\mathcal W}^{(n)} \bigr)$ and $\bigl({\mathcal W}, \tilde{{\mathcal W}}^{(n)} \bigr)$ are $(p(n),n)$-coupled webs. Observe that then
$\bigl({\mathcal W}^{(n)}, \tilde{{\mathcal W}}^{(n)} \bigr)$  forms a pair of $(\tilde{p}(n),n)$-coupled webs  where $1-p(n)= (1-\tilde{p}(n))^2$.

Consider $\Phi$  with ${\mathbf E} \bigl[ \Phi({\mathcal W})^2\bigr]<\infty$. Applying Theorem \ref{web:con} to $\bigl({\mathcal W}^{(n)}, \tilde{{\mathcal W}}^{(n)} \bigr)$ we obtain 
\begin{equation}
\label{strong1}
{\mathbf E} \Bigl[ {\mathbf E} \bigl[ \Phi({\mathcal W}^{(n)}) |{\mathcal W} \bigr]^2 \Bigr]={\mathbf E} \bigl[ \Phi({\mathcal W}^{(n)})\Phi(\tilde{{\mathcal W}}^{(n)})  \bigr]  \rightarrow {\mathbf E} \bigl[ \Phi({\mathcal W}^\prime)\Phi({\mathcal W}^{\prime\prime})  \bigr]
\end{equation}
where
$\bigl({\mathcal W}^\prime,{\mathcal W}^{\prime\prime}\bigr)$ is a $2\theta$-coupled pair of webs. 
Now by Theorem \ref{markov} we may assume that ${\mathcal W}, {\mathcal W}^\prime$, and  ${\mathcal W}^{\prime\prime} $ are defined on a common probability space with $\bigl({\mathcal W},{\mathcal W}^\prime\bigr)$ and $\bigl({\mathcal W},{\mathcal W}^{\prime\prime}\bigr)$ being $\theta$-couplings, and  ${\mathcal W}^\prime$ and ${\mathcal W}^{\prime\prime}$  conditionally independent given ${\mathcal W}$. Then 
\begin{equation}
\label{strong2}
{\mathbf E} \bigl[ \Phi({\mathcal W}^\prime)\Phi({\mathcal W}^{\prime\prime})  \bigr]= {\mathbf E} \Bigl[ {\mathbf E} \bigl[ \Phi({\mathcal W}^{\prime}) |{\mathcal W} \bigr]^2 \Bigr].
\end{equation}
Finally we have from the convergence of $\bigl({\mathcal W}, {\mathcal W}^{(n)}\bigr)$ to $\bigl({\mathcal W},{\mathcal W}^\prime\bigr)$ that for $\Psi({\mathcal W})$ satisfying ${\mathbf E} \bigl[ \Psi({\mathcal W})^2\bigr]<\infty$,
\begin{multline*}
{\mathbf E} \Bigl[ \Psi({\mathcal W}) {\mathbf E} \bigl[ \Phi({\mathcal W}^{(n)}) |{\mathcal W} \bigr] \Bigr]= \\
{\mathbf E} \bigl[ \Psi({\mathcal W}) \Phi({\mathcal W}^{(n)})  \bigr] \rightarrow {\mathbf E} \bigl[ \Psi({\mathcal W}) \Phi({\mathcal W}^{\prime})  \bigr]={\mathbf E} \Bigl[ \Psi({\mathcal W}) {\mathbf E} \bigl[ \Phi({\mathcal W}^{\prime}) |{\mathcal W} \bigr] \Bigr]
\end{multline*}
Thus taking $\Psi({\mathcal W})= {\mathbf E}\bigl[ \Phi({\mathcal W}^{\prime}) |{\mathcal W} \bigr]$
and  using \eqref{strong1} and \eqref{strong2} we obtain,
\begin{multline*}
{\mathbf E} \Bigl[ \bigl\{ {\mathbf E} \bigl[ \Phi({\mathcal W}^{(n)}) |{\mathcal W} \bigr]-{\mathbf E} \bigl[ \Phi({\mathcal W}^\prime) |{\mathcal W} \bigr]\}^2 \Bigr]= \\
{\mathbf E} \Bigl[ {\mathbf E} \bigl[ \Phi({\mathcal W}^{(n)}) |{\mathcal W} \bigr]^2 \Bigr] -2{\mathbf E} \Bigl[ \Psi({\mathcal W}) {\mathbf E} \bigl[ \Phi({\mathcal W}^{(n)}) |{\mathcal W} \bigr] \Bigr] +{\mathbf E} \Bigl[ {\mathbf E} \bigl[ \Phi({\mathcal W}^\prime) |{\mathcal W} \bigr]^2 \Bigr] \rightarrow 0.
\end{multline*}
\end{proof}

\section{The Erosion flow}

In this section we use the method  of filtering introduced by Le Jan and Raimond in Section 3 of \cite{lejan} to construct a stochastic flow of kernels from a pair of $\theta$-coupled Brownian webs. Some other examples of flows of kernels arising from filtering are  presented by  Le Jan and Raimond in \cite{lejan2}.

Let $ \bigl({\mathcal W}, {\mathcal W}^\prime \bigr)$ be a pair of $\theta$-coupled webs defined on a probability space $\bigl(\Omega, {\mathcal F}, {\mathbf P}\bigr)$. Writing $W(s,t,x)$ for $W_{x,s}(t)$, recall that as was observed in Proposition \ref{prop:ctsprod}, for each  $x\in{\mathbf R}$ and $s \leq t \leq u$, almost surely,
\begin{equation}
\label{mapcomp}
{ W}(t,u, { W}(s,t,x)) = { W}(s,u ,x).
\end{equation}
Now for each  $x\in{\mathbf R}$ and $s \leq t$, define $K_{s,t}(x,dy)$ to be (a version of) the conditional distribution of ${ W}_{x,s}(t)$ given ${\mathcal W}^\prime$, thus for all Borel subset $A \subseteq {\mathbf R}$,
\begin{equation}
\label{filter}
K_{s,t}(x,A) ={\mathbf P} \bigl( { W}_{x,s}(t) \in A |{\mathcal W}^\prime \bigr).
\end{equation}
We may choose a modification so that for each $s \leq t$ and Borel $A$, the map $(\omega,x) \mapsto K_{s,t}(\omega,x, A)$ is measurable, see Lemma 3.2 of \cite{lejan}.
Then using Fubini and Proposition \ref{joining}, we may deduce that  \eqref{mapcomp} implies that for each  $x\in{\mathbf R}$, $s \leq t \leq u$, almost surely for all Borel $A$,
\begin{equation}
 \int K_{s,t}(x,dy) K_{t,u} (y,A)= K_{s,u}(x,A).
\end{equation}
It is also easy to see that the family of random kernels  $\bigl(K_{s,t}; s\leq t \bigr)$ 
has stationary and  independent innovations, and  satisfies the definition of a measurable stochastic flow of kernels given in \cite{lejan}.

Our main aim in this section is to identify  the $N$-point motions of the flow of kernels just constructed by filtering. In general, Le Jan and Raimond have shown that the law of a stochastic  flow of kernels is characterized   by  its family of  $N$-point motions. For each integer  $N \geq 1$, the $N$ point 
motion 
of the flow $K$  is a Markov process on ${\mathbf R}^N$ with transition semigroup  given by
\begin{equation}
P^N_t(x,A)= {\mathbf E}\bigl[ K_{0,t}(x_1,A_1)K_{0,t}(x_2,A_2) \ldots K_{0,t}(x_N,A_N) \bigr],
\end{equation}  
for  $x=(x_1,x_2, \ldots, x_N) \in {\mathbf R}^N$ and $A=A_1\times A_2\times \cdots \times A_N $ a Borel cylinder set in  ${\mathbf R}^N$. 
Taking $N=1$ and the defining relation \eqref{filter}, we have
\begin{equation}
P^1_t(x,a)=  {\mathbf E}\bigl[ K_{0,t}(x,A)\bigr]= {\mathbf P} \bigl( { W}(0,t,x) \in A \bigr), 
\end{equation}
and hence the one point motion of $K$  is standard Brownian on ${\mathbf R}$.
The two-point motion can be identified with the help of the Markov property of $\theta$-couplings. Let $\bigl({\mathcal W},{\mathcal W}^\prime,{\mathcal W}^{\prime\prime}\bigr)$ be a triple of Brownian webs defined on a common probability space with each of $\bigl({\mathcal W}, {\mathcal W}^\prime\bigr)$ and $\bigl({\mathcal W}, {\mathcal W}^{\prime\prime}\bigr)$ being $\theta$-coupled Brownian motions, with ${\mathcal W}^\prime$ and ${\mathcal W}^{\prime\prime}$ being conditionally independent given ${\mathcal W}$. 
Then using this conditional independence,
\begin{multline}
P^2_t((x_1,x_2),A_1 \times A_2)=  {\mathbf E}\bigl[ K_{0,t}(x_1,A_1)K_{0,t}(x_2,A_2)\bigr] \\
={\mathbf E} \bigl[ {\mathbf P}\bigl({ W}^\prime(0,t,x_1)\in A_1 |{\mathcal W}\bigr) {\mathbf P}\bigl({ W}^{\prime\prime}(0,t,x_2)\in A_2 |{\mathcal W}\bigr) \bigr] \\
={\mathbf P}\bigl( W^\prime(0,t,x_1)\in A_1  \text{ and } { W}^{\prime\prime}(0,t,x_2)\in A_2 \bigr)= {\mathbf P}\bigl(X^\prime(t) \in A_1 \text{ and } X^{\prime\prime}(t) \in A_2 \bigr),
\end{multline}
where $\bigl(X^\prime, X^{\prime\prime} \bigr)$ are a pair of $2\theta$-coupled Brownian motions starting from $(x_1,x_2)$.

We now recall the main result from \cite{hw} concerning the characterization of consistent families of Brownian motions. We want to specify a diffusion in ${\mathbf R}^N$, where each  co-ordinate evolves as  Brownian motion, and each pair of co-ordinates forms a $\theta$-coupled pair of Brownian motions.   The  process may enter into parts of the state space where three or more co-ordinates are equal, and we need a way of describing the behaviour of the process when this occurs. The action of the (extended) generator of the process acting on $C^2$ does not characterize it; instead  we consider the action of the generator on certain space of piecewise linear functions. 

We begin by  partitioning  ${\mathbf R}^N$ into cells. A cell $E \subset {\mathbf R}^N$ is determined by some weak total ordering 
$\preceq$ of the $\{1,2,\ldots N\}$ via 
\begin{equation}
E= \{ x \in {\mathbf R}^N: x_i \leq x_j \text{ if and only if } i \preceq j \}.
\end{equation}
Thus $\{x \in {\mathbf R}^3: x_1 =x_2=x_3\}$, $\{x \in {\mathbf R}^3: x_1 <x_2=x_3\}$ and $\{ x\in {\mathbf R}^3: x_1>x_2>x_3\}$ are 
three of the thirteen distinct cells into which ${\mathbf R}^3$ is partitioned.

Suppose that $I$ and $J$ are disjoint subsets of $\{1,2,\ldots ,N\}$ with not both $I$ and $J$ empty. With such a pair we associate a 
vector $v=v_{IJ}$ belonging to ${\mathbf R}^N$ with components given by
\begin{equation}
\label{vdef}
v_i= \begin{cases} 0 & \text{ if $i \not\in I\cup J$,} \\
+1 & \text{ if $i \in I$,} \\
-1 & \text{ if $i \in J$.}
\end{cases}
\end{equation}
We want to associate with each point $x\in {\mathbf R}^N$ certain vectors of this form. To this end note that each point $x \in {\mathbf 
R}^N$ determines a partition $\pi(x)$ of $\{1,2,\ldots N\}$ such that $i$ and $j$ belong to the same component of 
$\pi(x)$ if and only if $x_i=x_j$. Then to each  point  $x \in {\mathbf R}^N$ we associate the set of vectors, denoted by ${\mathcal V}(x)$, which 
consists of
every vector of the form $v=v_{IJ}$ where $I\cup J$ forms  one component of the partition $\pi(x)$.

Let $ L_N$ be the space of real-valued functions defined on ${\mathbf R}^N$ which are continuous, and whose restriction to each cell is 
given by a linear function. Given a set of parameters $\bigl(\theta(k:l); k,l \geq 0\bigr)$ we define the operator ${\mathcal A}^\theta_N$ from $L_N$ to the space of 
real valued functions on ${\mathbf R}^N$ which are constant on each cell by
\begin{equation}
{\mathcal A}^\theta_N f (x) = \sum_{ v \in {\mathcal V}(x)} \theta(v) \nabla_{v} f(x).
\end{equation}
Here on the righthandside $\theta(v)= \theta(k:l)$ where $k=|I|$ is the number of elements in $I$ and $l=|J|$ is the number of elements in $J$ for $I$ and $J$  determined by $v=v_{IJ}$.  $\nabla_v f(x)$ denotes  the (one-sided) gradient of $f$ in the direction $v$ at the point $x$, that is 
\begin{equation}
\nabla_v f(x)= \lim_{\epsilon \downarrow 0}  \frac{1}{\epsilon} \bigl( f( x+\epsilon v)- f(x) \bigr).
\end{equation}

We say an ${\mathbf R}^N$-valued stochastic process  $\bigl( X(t); t \geq 0 \bigr)$ solves the ${\mathcal A}^\theta_N$-martingale problem if 
for  each $f \in L_N$,
\[
f\bigl(X(t)\bigr)- \int_0^t {\mathcal A}_N^\theta f \bigl(X(s)\bigr) ds \text{ is a martingale,}
\]
relative to some common filtration, and the bracket between co-ordinates $X_i$ and $X_j$ is given by
\[
\langle X_i,X_j\rangle (t)=\int_0^t {\mathbf 1}{( X_i(s)=X_j(s))} ds  \qquad \text{ for }  t\geq 0.
\]
In particular $\langle X_i \rangle(t)=t$.

The family of $N$-point motions of a stochastic flow of kernels will be consistent, in that, any $M \leq N$ co-ordinates taken from  the $N$-dimensional process  evolve as the $M$-dimensional member of the family. This translates into the following consistency condition on  $\bigl(\theta(k:l); k,l\geq 0 \bigr)$.  For any $k,l \geq 0$,
\begin{equation}
\label{eqnconsistent}
\theta(k:l)=\theta(k+1:l)+\theta(k:l+1), 
\end{equation}
Thinking of the $N$-dimensional process as describing the evolution of $N$-particles in ${\mathbf R}$,
 the co-efficient $\theta(k,l)$, for $k,l \geq 1$ can be thought of as describing the rate  at which  a cluster of  $k+l$ particles, momentarily moving together, separates into two clusters, one of $k$ particles and the other of $l$ particles.  
 Thus we expect the constraint, 
\begin{equation}
\label{eqnpositive}
\theta(k:l) \geq 0 \qquad \text{ for $k,l \geq 1$.}
\end{equation}  
In general this splitting of a cluster into two may  be asymmetric with $\theta(k:l) \neq \theta(l:k)$. In this case  the sizes of the two  clusters arising from a splitting  have an influence their relative ordering on the real line.   The parameters $\theta(k:0)$ and $\theta(l:0)$ introduce additional drift terms that compensate for such asymmetric splitting. In the  case of symmetric splitting, $\theta(k:l) = \theta(k:l)$ for all $k,l \geq 1$, and if there is no no-overall drift, $\theta(1:0)-\theta(0:1)=0$,  then the additional terms in the generator ${\mathcal A}_N^\theta$ corresponding to the $\theta(k:0)$ and $\theta(0:k)$ cancel each other and can thus be neglected entirely.

The main result  of \cite{hw} is the following.

\begin{theorem}
\label{main}
Let $\theta$ be a family of parameters satisfying the consistency and positivity properties given above in \eqref{eqnconsistent} and 
\eqref{eqnpositive}. For each $N \geq 1$ and $x \in {\mathbf R}^N$ there exists a process solving the ${\mathcal A}^\theta_N$-martingale 
problem 
starting from $x$. Moreover the law of this process is unique.
\end{theorem}

We are now able to identify the  $N$-point motions of the flow of kernels constructed  from a pair of $\theta$-coupled webs.
\begin{theorem}
\label{thm:motions}
The $N$-point motion of the stochastic flow of kernels derived from a pair of $\theta$-coupled webs solves the ${\mathcal A}_N^\theta$-martingale problem 
with $\bigl(\theta(k:l); k,l \geq 0 \bigr)$ given by
\begin{equation*}
\theta(k:l)= \begin{cases} 
\theta &\text{ if $k=1$ and $l=1$}, \\
\theta/2 &\text{ if either $k=1$ and $l \geq 2$ or vice versa,} \\
-m\theta/2 &\text{ if either $k=m+1\geq 2$ and $l=0$ or vice versa,}\\
0 & \text{ otherwise}.
\end{cases}
\end{equation*}
\end{theorem}

In preparation for proving this theorem we
let ${\mathbf P}^{\{k\}}_x$, for $x \in {\mathbf R}^N$ and $k \in \{1,2,\ldots N\}$, be probability measures governing an ${\mathbf R}^N$-valued  process $\bigl(Z(t); t \geq 0\bigr)$ such that under ${\mathbf P}^{\{k\}}_x$, $Z(0)=x$, and for all $i$ and $j$ not equal to $k$ the $i$th and $j$th components of $Z$ behave as a pair of coalescing Brownian motions, whilst the $k$th component is a Brownian motion independent of the others. Fix $f \in L_N$ and define real-valued functions $\psi_tf$ on ${\mathbf R}^N$ via
\begin{equation}
\label{defn:psi}
\psi^k_tf(x)= {\mathbf E}_x^{\{k\}} \bigl [ f(Z(t))- f(x)\bigr].
\end{equation}
and $\psi_t f= \sum_{k=1}^N \psi_t^k f$.
\begin{lemma}
\label{lem:dec}
Suppose $\theta(k:l); k,l \geq 0 \bigr)$ are as in the statement of the preceding theorem. Let $f \in L_N$ with
 ${\mathcal A}^\theta_N f  \geq 0$ everywhere. Then as $t \downarrow 0$, 
\[
\theta\sqrt{\frac{\pi}{t}}\psi_tf \downarrow {\mathcal A}^\theta_N f.
\]
\end{lemma}
\begin{proof}
Define closed  domains $D$ and $R_k$ for $k \in \{1,2,\ldots , N\}$ via
\begin{align*}
D&=\{ x \in {\mathbf R}^N: x_i\geq x_j \text{ for all } i \leq j \}, \\
R_k &= \{ x \in {\mathbf R}^N: x_i \geq x_j \text{ for all } i \leq j \text{ with } i,j \neq k \}.
\end{align*}
Then $ f \in L_N$ agrees on $R_k$ with a function $f_k$ of the form
\[
f_k(x)= \sum_{i \neq k} \beta^k_i|x_i-x_k| + \sum_i \alpha^k_i x_i.
\]
Recall we defined $\kappa_t(x)= {\mathbf E}\bigl[ |B(2t)+x|\bigr]- |x|$ where $B$ is a standard Brownian motion starting from zero. Then observing  that $R_k$ is an absorbing set for $X$ under ${\mathbf P}^{\{k\}}_x$, we obtain the representation 
\begin{equation}
\label{psirep}
\psi^k_tf(x)=\sum_{i \neq k} \beta^k_i\kappa_t(x_i-x_k) \qquad \text{ for } x \in D.
\end{equation}
Now consider ${\mathcal A}^\theta_Nf(x)$ for some $x \in D$. Suppose that $v \in {\mathcal V}(x)$ and that $v=v_{IJ}$  with either  $ I$ or $J$ equal to the singleton  $\{k\}$. Then  using  the representation of $f$ on $R_k$ we  obtain $\nabla_v f(x)+ \nabla_{-v} f(x)= 4\sum_{ i \neq k }\beta^k_i{\mathbf 1}(x_i=x_k)$. From this  we compute ${\mathcal A}^\theta_Nf(x)$ for the  values the parameters $\theta(k:l)$ given in the statement of the theorem, and we obtain,
\[
{\mathcal A}^\theta_Nf(x)=2 \theta \sum_{k, i \neq k } \beta_i^k {\mathbf 1}(x_i = x_k).
\]
Notice that the hypothesis that ${\mathcal A}^\theta_Nf \geq 0$ everywhere thus ensures that the co-efficients  $\beta_i^k$ are non-negative.

Now using \eqref{psirep} and the asymptotics for $\kappa_t(x)$ given at \eqref{kappalimit} we obtain 
\[
\theta\sqrt{\frac{\pi}{t}}\sum_k \psi^k_tf (x) \downarrow 2\theta \sum_{k, i \neq k}\beta_i^k {\mathbf 1}(x_i = x_k)={\mathcal A}^\theta_Nf(x).
\]
Finally, note that, by applying a suitable permutation to the co-ordinates,  the results extends from all $x \in D$, to all $x \in {\mathbf R}^N$.
\end{proof}

\begin{proof}[Proof of Theorem \ref{thm:motions}]
Let ${\mathcal W}$, ${\mathcal W}^\prime$ and ${\mathcal W}^{(n)}$, for $n \geq 1$ be Brownian webs defined on a common probability space and such that
$\bigl({\mathcal W}, {\mathcal W}^{(n)}\bigr)$ are $(p(n), n)$-coupled where $p(n)= \theta \sqrt{\pi/n}$, and $\bigl({\mathcal W}, {\mathcal W}^{\prime}\bigr)$ are $\theta$-coupled.
By  Proposition \ref{prop:strong} we have
\[
{\mathbf P} \bigl( { W}^{(n)}(0,t,x) \in A | {\mathcal W} \bigr) \stackrel{L^2}{\rightarrow} {\mathbf P} \bigl( { W}^{\prime}(0,t,x) \in A | {\mathcal W} \bigr),
\]
Since the lefthandside is uniformly bounded in $L^\infty$, it follows that,
\[
\prod_{k=1}^N {\mathbf P} \bigl( { W}^{(n)}(0,t,x_k) \in A_k | {\mathcal W} \bigr) \stackrel{L^2}{\rightarrow} \prod_{k=1}^N {\mathbf P} \bigl( { W}^{\prime}(0,t,x_k) \in A_k | {\mathcal W} \bigr).
\]
Thus the semigoup of the $N$-point motion of the  erosion flow is given by
\begin{multline*}
 {\mathbf E}\left[ \prod_{k=1}^N K_{0,t}(x_k,A_k) \right] =\lim_{n \rightarrow \infty} {\mathbf E} \left[\prod_{k=1}^N {\mathbf P} \bigl( { W}^{(n)}(0,t,x_k) \in A_k | {\mathcal W} \bigr) \right] \\
 =\lim_{n \rightarrow \infty} {\mathbf P} \bigl( { W}^{(k,n)}(0,t,x_k) \in A_k  \text{ for } k=1,2,\ldots,N \bigr),
\end{multline*}
where for each $n$, ${\mathcal W}^{(k,n)}$ for  $k=1,2,\ldots,N$ are Brownian webs, conditionally independent given ${\mathcal W}$, and  such that $\bigl( {\mathcal W}^{(k,n)}, {\mathcal W} \bigr)$ are $(n, p(n))$-coupled. Consequently we see that we must prove that as $n$ tends to infinity the ${\mathbf R}^N$-valued process 
\begin{equation}
\label{motionsinwebs}
\bigl( {\mathcal W}^{(1,n)}(0,x_1,t), {\mathcal W}^{(2,n)}(0,x_2,t), \ldots, {\mathcal W}^{(N,n)}(0,x_N,t); t \geq 0 \bigr)
\end{equation}
 converges in distribution to the unique solution to the ${\mathcal A}^\theta_N$-martingale problem starting from $(x_1, x_2, \ldots, x_N)$. 
 
We may construct an ${\mathbf R}^N$-valued process having the same law as that at \eqref{motionsinwebs} by generalizing the construction of $(p,n)$-coupled Brownian motions.  Let $ \bigl(Y_k^i;   k \geq 0, 1 \leq i \leq N \bigr)$ be a family of independent Bernoulli($p(n)$) random variables. Let $S_k=\{ i : Y_k^i=1\}$. Conditionally on $ \bigl(Y_k^i;   k \geq 0, 1 \leq i \leq N \bigr)$ let $\bigl(X^{(n)}_1(t), X^{(n)}_2(t),\ldots X({(n)}_N(t); t \geq 0\bigr)$ be Brownian motions starting from $(x_1,x_2 \ldots ,x_N)$
with $X^{(n)}_i$ evolving independently of $\bigl(X^{(n)}_j; j \neq i \bigr)$ during the time  interval $t \in [k/n,(k+1)/n]$ if $i \in S_k$,
but with $X_i^{(n)}$ and $X_j^{(n)}$ coalescing if they meet during the interval $t \in [k/n,(k+1)/n]$ and both $i,j \not\in S_k$.
 
 The  sequence laws of  the processes $\bigl(X^{(n)}; n \geq 1 \bigr)$ is tight since each $X^{(n)}_k$ is a Brownian motion. Thus it is sufficient to show any subsequence $\bigl(X^{(n_k)}; k \geq 1 \bigr)$ which converges in distribution must converge to the unique solution of the ${\mathcal A}^\theta_N$-martingale problem.  We may assume for notational simplicity  that  $\bigl(X^{(n)}; n \geq 1\bigr)$ itself converges in distribution to a process $X$.
Each co-ordinate $X_k$ of $X$ must be distributed as a Brownian motion, and moreover be a martingales relative to the natural filtration of $X$, since the corresponding statement holds for  $X^{(n)}$. Furthermore each pair of co-ordinates $\bigl(X_i,X_j\bigr)$ must be a $2\theta$-coupled Brownian motions, by Proposition \ref{bm:con1}. It remains to show that the process $f\bigl(X(t)\bigr)- \int_0^t {\mathcal A}^\theta_Nf\bigl(X(s)\bigr)ds$ is a martingale for each $f \in L_N$.

Fix $f \in L_N$ satisfying ${\mathcal A}^\theta_Nf \geq 0$ everywhere. Let $\bigl({\mathcal F}^{(n)}_t; t \geq 0 \bigr)$ denote the natural  filtration of the process $\bigl(X^{(n)}(t); t \geq 0 \bigr)$. We define  $h_n(x)$ via
\[
h_n(x)=\sum_{S \subseteq \{1,2,\ldots , n\}} p(n)^{|S|} (1-p(n))^{N-|S|}{\mathbf E}^S_x \bigl[ f(Z({1/n}))-f(x)\bigr],
\]
where the process $Z$, governed by probability measure ${\mathbf P}^S_x$, starts from $x$ and evolves so thats  its $i$th and $j$th components with  $i,j \not\in S$  are coalescing Brownian motions, whilst its $i$th component for $i \in S$ is  a Brownian motion independent of the other components. We observe that it follows easily from the construction of $X^{(n)}$ just given that the discrete parameter process
\begin{equation}
\label{martingale}
f\bigl(X^{(n)}({k/n})\bigr) - \sum_{r=0}^{k-1} h_n\bigl((X^{(n)}(r/n)\bigr)
\end{equation}
is a martingale, relative to the discrete parameter filtration $ \bigl( {\mathcal F}^{(n)}_{k/n}; k\geq 0 \bigr)$. 
Using the fact that $f$ is Lipshitz we have
\[
\bigl| {\mathbf E}^S_x \bigl[ f(Z({1/n}))-f(x)\bigr] \bigr| \leq \frac{C}{\sqrt{n}},
\]
for a suitable constant $C< \infty$ depending on $f$ only. Then
\[
h_n(x) \leq \sum_{k=1}^N {p(n)}  {\mathbf E}^{\{k\}}_x \bigl[ f(Z({1/n}))-f(x)\bigr] +\frac{C}{\sqrt{n}} N^2p(n)^2.
\]
Recalling that $p(n)=\theta \sqrt{\pi/n}$ and the definition of $\psi_tf$, we see that for any $\epsilon>0$, and   all sufficiently large $n$, 
\[
h_n(x) \leq \theta \sqrt{\frac{\pi}{n}}\psi_{1/n}f(x) + \frac{\epsilon}{n} \leq \frac{1}{n} \left\{ \theta (\pi m)^{1/2}\psi_{1/m}f(x) + {\epsilon}{}\right\},
\]
the second inequality following because, by Lemma \ref{lem:dec}, $ t \mapsto t^{-1/2} \psi_t f(x)$ decreases as $t$ does.
Returning to the martingale given at \eqref{martingale}, whenever the above inequality holds,
\[
f\bigl(X^{(n)}({k/n})\bigr) -\frac{1}{n} \sum_{r=0}^{k-1} \left\{ \theta (\pi m)^{1/2}\psi_{1/m}f\bigl(X^{(n)}(r/n)\bigr) + \epsilon\right\}
\]
is a supermartingale, from which it follows by standard weak convergence arguments similar to in Proposition \ref{bm:con1} that
\[
f\bigl(X^{}({t})\bigr) -\int_0^t \left\{ \theta (\pi m)^{1/2}\psi_{1/m}f\bigl(X^{}(s)\bigr) + \epsilon\right\}ds
\]
is a supermartingale relative to the natural filtration of $X$. Now letting $\epsilon$ tend to zero, and $m$ tend to infinity, appealing once again to Lemma \ref{lem:dec}, we deduce that
$f\bigl(X(t)\bigr)- \int_0^t {\mathcal A}^\theta_Nf\bigl(X(s)\bigr)ds$ is a supermartingale also.

To complete the prove of the theorem we consider a general $f \in L_N$. Put 
\[
g(x)=  \sum_{i<j} | x_i- x_j|.
\]
Then for sufficiently large $k>0$ we have both ${\mathcal A}^\theta_N (kg+f) \geq  0$ and ${\mathcal A}^\theta_N(kg-f) \geq 0 $, and thus we may apply the previous result to both $kg+f$ and $kg-f$. However we know that for each choice of distinct $i,j$ that the process $ \bigl( X^i(t), X^j(t); t \geq 0 \bigr)$ is a 
$2\theta$-joining of Brownian motions, and consequently
\[
g\bigl(X(t)\bigr)- 4 \theta \sum_{i < j} \int_0^t {\mathbf 1}\bigl(X_i(s)=X_j(s)\bigr) ds 
\]
is a martingale. Observing that 
\[
{\mathcal A}^\theta_N g (x) = 4\theta \sum_{i < j} {\mathbf 1}(x_i=x_j),
\]
we deduce that  $f\bigl(X(t)\bigr)- \int_0^t {\mathcal A}_N^\theta f\bigl(X(s)\bigr)ds$ is both a supermartingale and a submartingale.
\end{proof}

\section{Generalisations}

Whilst we were writing this paper we became aware of the work of  Sun and Swart \cite{ss} on the Brownian net. A key ingredient in their construction of the net is a pair of coupled Brownian webs which they call left/right webs. In fact $\theta$-coupled webs and left/right webs are both special cases of a three parameter family of couplings. The following proposition generalizes both Proposition \ref{bm:unique} and Proposition 16 of \cite{ss}. 
It is important  note that in this result $\tilde{L}^0_t(Z)$ denotes the symmetric local time of the semimartingale $Z$ at zero, thus
$ \tilde{L}^0_t(Z)=\tfrac{1}{2}\bigl(L^0_t(Z) +L^{0-}_t(Z)\bigr)$ where $L^{x}_t(Z)$ denotes the usual right continuous version of local times of $Z$.

\begin{proposition}
\label{prop:gencouple}
Suppose that $\beta_1, \beta_2$ and $\theta$ are parameters satisfying $|\beta_1-\beta_2|\leq 2 \theta<\infty$. Then, for each given starting point $ (x_1,x_2) \in {\mathbf R}^2$, there exists a stochastic process $ \bigl( (X(t),X^\prime(t); t \geq 0 \bigr)$ such that $X$ is a Brownian motion with drift $\beta_1$ starting from $x_1$, and $X^\prime$ is a Brownian motion with drift $\beta_2$ starting from $x_2$ (relative to some common filtration), and 
\begin{align}
\langle X,X^\prime\rangle (t)=\int_0^t {\mathbf 1}_{( X(s)=X^\prime(s))} ds  \qquad &  t\geq 0, \\
\label{local}
\tilde{L}^0_t\bigl( X-X^\prime\bigr) = 2\theta \int_0^t {\mathbf 1}_{( X(s)=X^\prime(s))} ds  \qquad& t\geq 0.
\end{align}
Moreover the law of $\bigl(X,X^\prime \bigr)$ is uniquely determined.
\end{proposition}

\begin{proof}
In  case $\beta_1=\beta_2=\theta=0$, $\bigl( X,X^\prime \bigr)$ are a pair of coalescing Brownian motions. Henceforth we assume $\theta>0$.

For $x \in {\mathbf R}$, $d \in {\mathbf R}$, and $\beta \in [-1,1]$, the stochastic differential equation
\begin{equation}
\label{sde}
Z(t)= x+B(t) + dt+ \beta \tilde{L}^0_t(Z), \qquad t \geq 0,
\end{equation}
where  $Z$ is assumed to be adapted to a filtration relative to which $B$ is a  Brownian motion,
possesses a pathwise unique solution. In particular   the solution is  unique-in-law. The case  $d=0$ is known as skew Brownian motion with skewness parameter $\beta$ and general $d$ can be reduced to this by a change of measure. 

 To construct $\bigl(X,X^\prime\bigr)$ with  the  properties stated in the proposition, we let $Z$ be a solution to the preceding SDE with parameters given by $x=(x_1-x_2)/2$, $d=(\beta_1-\beta_2)$ and $\beta=(\beta_1-\beta_2)/(2\theta)$, and let $B^\prime$ be an independent Brownian motion defined on the same probability space as $Z$. We define a continuous strictly increasing process  by
\[
 \alpha(t)= {2}t + \frac{1}{\theta}\tilde{L}^0_t(Z),
\]
 and let $\bigl(A(t);  t \geq 0 \bigr)$ be the inverse of $ \bigl( \alpha(t); t \geq 0 \bigr)$. 
 We define $X$ and $ X^\prime$ via, 
\begin{align}
\label{xfromb}
X(t) &=x_1+ B^\prime\bigl( t-A(t) \bigr)+B\bigl(A(t)\bigr) +\beta_1 t \\
\label{xfromb2}
X^\prime(t) &=x_2+ B^\prime\bigl( t-A(t) \bigr)-B\bigl(A(t)\bigr) +\beta_2 t.
\end{align}
Then $X$ and $X^\prime$ are Brownian motions with drifts $\beta_1$ and $\beta_2$ relative to the filtration they jointly generate. Moreover
\[
\langle X, X^\prime \rangle (t)=t -2A(t)=\frac{1}{\theta}\tilde{L}^0_{A(t)}(Z),
\]
the second equality holding by virtue of   the definition of  $\alpha$. 
Now computing $X(t)- X^\prime(t)$ we obtain,
\[
X(t)- X^\prime(t)=(x_1-x_2)+2 B\bigl(A(t)\bigr)+ (\beta_1-\beta_2)t= 2Z\bigl(A(t)\bigr),
\]
which implies on the one hand that that $\tilde{L}^0_t(X-X^\prime)=2\tilde{L}^0_{A(t)}(Z)$, and on the other that
\[
\int_0^t {\mathbf 1}\bigl(X(s)=X^\prime(s) \bigr) ds = \int_0^{A(t)}{\mathbf 1} \bigl(Z(s)=0\bigr) d\alpha(s)= \frac{1}{\theta}\tilde{L}^0_{A(t)}(Z).
\]
Combining the two gives the desired expression for  $\tilde{L}^0_t(X-X^\prime)$, and also  for $\langle X, X^\prime \rangle (t)$. Thus existence is proved.

Now assume that $\bigl(X,X^\prime\bigr)$ is any pair of processes having the properties specified  in the proposition.   
The processes defined by $U(t)= \bigl(X(t)+X^\prime(t)\bigr)/2 -(x_1+x_2)/2 - (\beta_1+\beta_2)t/2$ and $V(t)=\bigl(X(t)-X^\prime(t)\bigr)/2 -(x_1-x_2)/2 - (\beta_1-\beta_2)t/2$ are orthogonal martingales with  $\langle U \rangle (t)=t-C(t)$ and $\langle V \rangle  (t)= C(t)$ where
\[
C(t)=\frac{1}{2}\int_0^t {\mathbf 1}\bigl(X(s)\neq X^\prime(s) \bigr) ds.
\]
Consequently by Knight's theorem  we may represent $U$ and $V$ as time changes of independent Brownian motions: $U(t)= B^\prime\bigl(t-C(t)\bigr)$ and $V(t)= B\bigl(C(t) \bigr)$. Expressing $X$ and $X^\prime$ in terms of $U$ and $V$ now gives \eqref{xfromb} and \eqref{xfromb2}, but with $C(t)$ in place of $A(t)$.
We observe that $\bigl( C(t); t \geq 0 \bigr)$ is continuous and strictly increasing, for if it were constant on some interval then $X(t)=X^\prime(t)$ would hold on the interval, and this would not be consistent with the expression for the local time of $X-X^\prime$.  Let $\bigl(\gamma(t); t \in [0,C(\infty)) \bigr)$ be its inverse, and set $Z(t)=\bigl(X(\gamma(t))-X^\prime(\gamma(t))\bigr)/2$. Then $\tilde{L}^0_t(X-X^\prime)=2\tilde{L}^0_{C(t)}(Z)$ and we easily check that $Z$ satisfies \eqref{sde}. In particular $C(\infty)=\infty$ almost surely.
Moreover we may now identify $C$ with the process $A$ which appears in the proof of existence, and then \eqref{xfromb} and \eqref{xfromb2} hold for $\bigl(X,X^\prime\bigr)$.  Consequently the law of $\bigl(X,X^\prime\bigr)$ is uniquely determined.
\end{proof}

Following the same steps as in Section 3 of this paper, which are also described in \cite{ss}, we are able to construct a pair of coupled webs $\bigl({\mathcal W}, {\mathcal W}^\prime \bigr)$ such that each path in ${\mathcal W}$ is a Brownian motion drift $\beta_1$, each path in ${\mathcal W}^\prime$ is a Brownian motion drift $\beta_2$, and a pair of paths, one from each of the webs, evolves as the diffusion specified by the previous proposition. Then we may again define a flow of kernels from this pair of webs via \eqref{filter}.
We conjecture that the $N$-point motion of this flow  solves the ${\mathcal A}_N^\theta$-martingale problem 
with $\bigl(\theta(k:l); k,l \geq 0 \bigr)$ given by
\begin{equation*}
\theta(k:l)= \begin{cases} 
\theta &\text{ if $k=1$ and $l=1$}, \\
 ( 2\theta + \beta_1 - \beta_2)/4,& \text{ if $k=1$ and $l \geq 2$, }\\
(2 \theta  + \beta_2-\beta_1)/4,& \text{ if $k\geq 2$ and $l =1$, }\\
 0& \text{ if both $k,l \geq 2$},
\end{cases}
\end{equation*}
with $\theta(k:0)$ and $\theta(0:l)$ determined by the consistency rule together with $-\theta(0:1) = \theta(1:0) = \beta_1/2$. These are, in general, asymmetric  flows, in which we imagine  erosion to be occurring  at different rates on the right and left sides of a cluster.

\end{document}